\documentclass[preprint,11pt,authoryear]{elsarticle}
\usepackage{graphicx}
\usepackage[utf8]{inputenc}
\usepackage{amsmath}
\usepackage{algorithm, algorithmicx}
\usepackage[noend]{algpseudocode}
\algrenewcommand{\algorithmiccomment}[1]{\hfill// #1}
\algnewcommand{\LineComment}[1]{\State // #1}
\algnewcommand\algorithmicto{\textbf{\textup{ to }}}
\algnewcommand\Input{\item[\textbf{Input:}]}
\algnewcommand\Output{\item[\textbf{Output:}]}
\algnewcommand{\IfThenElse}[3]{% \IfThen{<if>}{<then>}
	\State \algorithmicif\ #1\ \algorithmicthen\ #2 \ \algorithmicelse\ #3
}
\algnewcommand\Execute{\State \textbf{execute} }
\algnewcommand\Initialize{\State \textbf{initialize} }

\usepackage{longtable}
\usepackage[inline, shortlabels]{enumitem}
\usepackage{multirow,booktabs,adjustbox,lipsum}
\usepackage{amssymb,mathtools,amsthm,thmtools,bbm}
\usepackage{eurosym}
\DeclareMathOperator*{\argmax}{argmax} % thin space, limits underneath in displays

\theoremstyle{definition}

\theoremstyle{remark}

\usepackage[margin=1in]{geometry}
\usepackage{graphicx}
\usepackage{caption}
\usepackage{subcaption}

\usepackage{soul, multirow, bm, amssymb, mathtools, array, enumitem}

\usepackage{tikz}
\usetikzlibrary{shapes,arrows,positioning}
\usetikzlibrary{automata,arrows,positioning,calc}
\usetikzlibrary{decorations.pathreplacing}

\usepackage{lipsum}
\tikzset{naming/.style={align=center,font=\footnotesize}}
\tikzset{area/.style = {draw, shape = regular polygon, regular polygon sides = 6, thick, minimum width = 5cm}}

\usepackage[colorlinks=true, urlcolor=blue]{hyperref}

\usepackage{xcolor}

\usepackage{natbib}
 \bibpunct[, ]{(}{)}{,}{a}{}{,}%
 \def\bibsep{\smallskipamount}%
\usepackage{eurosym}
\usepackage{lineno}

\journal{European Journal of Operational Research}

\usepackage{xcolor}

\begin{document}

\begin{frontmatter}



\title{Strategic Selection of Remanufacturing Business Models: A Consumer Perception Perspective}


\author[1]{Zhongxin Hu}
\ead{z.hu1@tue.nl}
\author[1]{Christina Imdahl}
\ead{c.imdahl@tue.nl}
\author[1]{Z\"{u}mb\"{u}l Atan\corref{cor1}}
\ead{z.atan@tue.nl}
\cortext[cor1]{Corresponding Author.}

\affiliation[1]{organization={ Department of Industrial Engineering and Innovation Sciences, Eindhoven University of Technology},%Department and Organization
            city={Eindhoven},
            country={The Netherlands}}

\begin{abstract}
As a key circular economy strategy, remanufacturing allows original equipment manufacturers (OEMs) to reduce waste by restoring used products to ``as-new'' conditions. This paper investigates an OEM's optimal remanufacturing business model by incorporating consumer perceptions into price and production quantity decisions. We analyze three alternative models: no remanufacturing, OEM in-house remanufacturing, and third-party remanufacturer (TPR) authorized remanufacturing. We extend the authorization with a two-part tariff contract and consider a stochastic market size. Through a numerical approach, we optimize price and quantity decisions based on consumer perceptions and develop a hierarchical decision roadmap to guide model selection. Our findings show that when consumer's perceived value of remanufactured products is high, OEM in-house remanufacturing is most profitable and reduces environmental impacts, but generally leads to a market dominated by remanufactured products. In contrast, when consumer's perceived value of remanufactured products is moderate and TPR remanufacturing significantly increases the perceived value of new products, the TPR-authorized remanufacturing is most profitable. It typically boosts total market sales, but accordingly increases environmental impacts. In addition, sensitivity analysis indicates that two-part authorization contracts are more advanced in meeting stringent environmental requirements than one-part contracts. Incorporating market size stochasticity enhances system profitability while keeping environmental impacts within a limited scope.
\end{abstract}


% \begin{highlights}
% \end{highlights}


\end{frontmatter}

\section{Introduction}\label{sec:Intro}
In recent years, remanufacturing has emerged as a key strategy within the circular economy, offering substantial potential for resource conservation, waste reduction, and the creation of economic value. Remanufacturing involves processes that restore used products to an ``as-new'' state, effectively prolonging their useful lives.

Environmentally, remanufacturing reduces the need for new raw materials, decreases energy consumption and emissions compared to producing brand-new products, and minimizes waste by diverting products from landfills back into the economic cycle. \citet{units4.1} shows that remanufacturing in the automotive sector can decrease virgin material use by 88\%, reduce energy requirements by 56\%, and cut CO$_2$ emissions by up to 53\%. Economically, it presents substantial advantages by potentially saving up to 40\% - 64\% of production costs compared to manufacturing new products \citep{ginsburg2001once,chatti2019cirp}. \citet{WorldBank2022} reports that the European remanufacturing market is currently valued at \euro31 billion and could grow to \euro100 billion by 2030, saving 21 megatons of CO$_2$ emissions.

Given these advantages, a variety of remanufacturing business models have been developed, including in-house remanufacturing by the original equipment manufacturer (OEM) \citep{lin2024house, abbey2024closed}, outsourcing to a third-party remanufacturer (TPR) \citep{wang2017design}, and authorizing to a TPR \citep{zou2016third, banerjee2023optimal}. In both OEM in-house and outsourcing models, the OEM maintains control over the sale and marketing of remanufactured products, but the authorizing model allows TPR to sell and market remanufactured products under its own brand. This distinction in branding improves the consumer's awareness of who is responsible for the remanufactured product. 

The effectiveness of the three business models depends not only on operational and economic factors, but also on consumer perception \citep{subramanian2012key,donohue2020behavioral,huang2024pride}. Empirical evidence shows that consumers generally value remanufactured products less than new products, with a perceived value ranging between 40\% and 90\%, the so-called willingness-to-pay (WTP) discount factor \citep{guide2010potential, abbey2017role}. \citet{agrawal2015remanufacturing} find the perceived value of new products can be significantly influenced by the identity of the remanufacturer. Specifically, for products, such as MP3 players and printers, third-party remanufacturing can increase the perceived value of new products by up to 8\%, the so-called \textit{contrast effect}, while OEM remanufacturing may decrease it by up to 7\%, the so-called \textit{assimilation effect}. Although some analytical models have partially incorporated consumer perceptions \citep{agrawal2010essays,fang2020third,li2024should}, most research treats these perceptions as fixed assumptions and primarily focuses on cost analysis, assessing how cost parameters determine business model choices given a specific condition of consumer perceptions. Only a limited number of studies have systematically explored the patterns by which consumer perception influences the optimal business model. To address this research gap, we investigate the following research question: \textit{How does consumer perception affect the choice of the optimal remanufacturing business model?} 

In this study, we compare three remanufacturing business models: no remanufacturing (Model N), OEM in-house remanufacturing (Model O), and TPR-authorized remanufacturing (Model T). The outsourcing model, in which the OEM retains control over sales and marketing, is equivalent to OEM in-house remanufacturing from a consumer perspective; thus, it is excluded from analysis. Using a numerical method, we determine unit prices and quantities of new and remanufactured products to be offered to consumers so that OEM's profit is maximized. Subsequently, we identify the most profitable model under varying consumer perception conditions and reveal systematic patterns in model selection. Further, we assess both market and environmental outcomes under the optimal model, and discuss the assumptions of authorization contracts and market size.

We contribute to the remanufacturing literature in multiple ways. First, we integrate consumer perceptions into remanufacturing business models and examine their impact on the selection of the best business model. We develop a hierarchical decision roadmap that accommodates the full spectrum of consumer perceptions.
Second, we examine the market and environmental outcomes associated with optimal business model selection, thereby evaluating the extent to which remanufacturing can promote sustainability. Third, we study stochastic market size and generalize the authorization contract by a two-part tariff structure, including a one-time fixed fee and a unit fee. 

Our findings reveal that consumer perceptions substantially affect the OEM's selection of the remanufacturing business model. Among perception factors, the WTP discount factor for remanufactured products exerts the strongest influence, whereas effects related to remanufacturer identity, i.e., assimilation and contrast effects, influence only under specific conditions. When consumers perceive remanufactured products as less valuable, that is, low WTP discount factors, the OEM should avoid remanufacturing. For moderate WTP discount factors, the optimal model depends on the magnitude of assimilation and contrast effects. When these effects are moderate to high, authorizing a third-party remanufacturer is the most profitable option. When consumers perceive remanufactured products as nearly equivalent to new ones, that is, high WTP discount factors, OEM in-house remanufacturing is most profitable. Based on these insights, we develop a hierarchical roadmap to guide the selection of the optimal remanufacturing business model aligned with varying consumer perception conditions.

Our numerical analysis explores the market dynamics and environmental impact of selected optimal models. 
When the TPR-authorized remanufacturing model becomes optimal, the total product output can increase by up to 63.9\% and new product output can rise by up to 14.8\%, alleviating managerial concerns about cannibalization \citep{guide2010potential,atasu2010so}. However, this substantial market expansion also increases overall consumption and emissions, outweighing the environmental benefits of remanufacturing and resulting in a higher aggregate environmental impact. 
Consequently, the TPR-authorized model generates a ``win-win-lose'' outcome: higher profits and significant market growth, but an increased environmental burden. In contrast, when the OEM in-house remanufacturing model is optimal, a ``win-win-win'' outcome emerges: higher profits, moderate market growth, and reduced environmental impact. In this case, remanufactured products can dominate the market when consumers perceive their value as sufficiently high, resulting in a market where only remanufactured goods are offered.

Through sensitivity analysis, we investigate how the design of authorization contracts and the stochasticity of market size affect system profitability and environmental outcomes. Numerical results across a realistic range of consumer perceptions suggest that although two-part tariff contracts do not always generate the highest profits, they offer greater flexibility than one-part contracts in satisfying stringent environmental requirements. Introducing stochastic market size not only increases profit by up to 44.8\%, but also reduces environmental impact in some cases; even when environmental impact increases, the increase typically remains below 5\%. 

The remainder of the paper is organized as follows. Section~\ref{sec:Literature review} reviews the relevant literature. Section~\ref{sec:model framework} introduces our model framework. Section~\ref{sec:remanufacturing modes} sets up three alternative models and presents their optimal and approximate solutions. Section~\ref{sec:selection} compares the equilibrium results across models using numerical experiments and analyzes their implications for market dynamics and environmental impact. Section~\ref{sec:contract and stochasticity} discusses the impact of the authorization contract and the stochasticity of the market size. Section~\ref{sec:conclusion} concludes with implications and points out further research directions.

\section{Literature review} \label{sec:Literature review}

Our research contributes to two streams of literature in operations management: remanufacturing business model and consumer perceptions in remanufacturing. 

\subsection{Remanufacturing business model}
Recent research on OEM's selection of remanufacturing business model generally distinguishes between two main categories: OEM in-house remanufacturing and TPR remanufacturing, which can be achieved through outsourcing or authorization. The OEM in-house remanufacturing model offers OEM control over brand image, process technology, and product acquisition, potentially leading to improved environmental impact \citep{orsdemir2014competitive} and operational efficiency \citep{lin2024house}. This model is exemplified by companies such as Xerox and Caterpillar, but requires significant technical and financial investment \citep{fang2020third, abbey2024closed}. In contrast, cooperation with a TPR allows the OEM to delegate remanufacturing activities, either by outsourcing or authorization. In outsourcing, the OEM retains the sale and marketing of remanufactured products by the original brand, while authorization enables third parties to obtain license and technical support and sell remanufactured products under their own brands. 
The distinctions in the responsibility for sales and marketing can affect the willingness of consumers to purchase remanufactured goods and, therefore, should be taken into account when comparing remanufacturing business models. 

Current comparative studies of remanufacturing business models emphasize mainly cost analysis, evaluating how remanufacturing costs impact the selection of the most profitable model \citep{zou2016third,feng2021environmentally,liu2022rent}. %However, the distinctions in sales and marketing and their impact on consumer perceptions and demand for new and remanufactured products remain less explored in the literature of remanufacturing models.
Most of the existing literature on third-party remanufacturing \citep{huang2024remanufacturing,li2024should} assumes a one-part contract with a unit fee. Motivated by the practice of licensing and patent agreements \citep{bonnet2006two,san2015optimal,banerjee2023optimal}, this study introduces a two-part structure to better analyze the design of contracts. The two-part contract includes a fixed one-time fee for technical permission and license rights, and a unit fee for each remanufactured product. 

\subsection{Consumer perception in remanufacturing}
Emerging research emphasizes the importance of considering customer interactions in business decision-making processes, particularly in circular economy practice \citep{subramanian2012key,abbey2017role,donohue2020behavioral,huang2024pride}. 

Empirical evidence shows that consumers generally perceive remanufactured products as less valuable than new products, with the value gap generally ranging from 40\% to 90\% depending on the product categories \citep{guide2010potential,subramanian2012key,abbey2015remanufactured,abbey2017role}.
This insight is incorporated into analytical models of remanufacturing in various ways. Most studies employ willingness-to-pay (WTP) to measure consumers' perceived value, typically assuming a constant WTP discount factor for remanufactured products relative to new ones \citep{orsdemir2014competitive,zou2016third, fang2020third, shi2020distribution, liu2022rent}. Only a few studies deviate from this convention by a variable WTP discount factor \citep{abbey2017role} or by a demand-switching fraction \citep{ovchinnikov2011revenue}.

Recent empirical research from \citet{agrawal2015remanufacturing} shows that the presence of remanufactured products and the identity of the remanufacturer can influence consumers' perceived value of new products. Specifically,  for consumer products, such as MP3 players and printers, if both new and remanufactured products are provided from the same manufacturer, i.e., the OEM, consumers tend to view their value as similar, which reduces the perceived value of new products by up to 7\% due to the \textit{assimilation effect}. In contrast, when remanufactured products are provided by a third-party remanufacturer, consumers see new and remanufactured products as different, increasing the perceived value of new products by up to 8\% through the \textit{contrast effect}. 

Several studies integrate assimilation and contrast effects in the selection of remanufacturing business models. For example, \cite{agrawal2010essays} analytically examines these effects in a context where the OEM collects used products preemptively and a competitive TPR remanufactures independently. 
Building on this, \citet{fang2020third} investigate the decision making of an OEM under a competitive third-party market entry and conclude that OEM involvement in remanufacturing is optimal only when both the remanufacturing costs and contrast effect are low. \citet{wu2020competitive} study a market with two competing OEMs with different brand equity and find that the presence of assimilation and contrast effects diminished the incentive of both firms to remanufacture. Although these studies incorporate consumer perceptions into remanufacturing business models, they do not cover scenarios where the OEM cooperates strategically with the TPR. To this end, \citet{huang2024remanufacturing} compare OEM self-operation remanufacturing and licensing with a TPR, but do not consider the effects of consumer perceptions. Likewise, \citet{li2024should} analyze the setting in which the OEM authorizes a TPR using a one-part contract and account for consumer perceptions, but do not directly compare the impacts of assimilation and contrast effects side by side, which could have provided deeper insights into dynamic trade-off involved in model selection. 

To address these gaps, we evaluate OEM in-house remanufacturing with assimilation effect and TPR-authorized remanufacturing with contrast effect. 
Most importantly, unlike prior studies that focused on cost analysis, we analyze how consumer perceptions can shift the optimal model selection. In addition, although the existing literature on remanufacturing business models has considered stochastic demand related to heterogeneous consumers, we introduce a stochastic market size, that is, allow the total number of potential consumers to vary. To our knowledge, this is the first paper to jointly incorporate customer heterogeneity and market size stochasticity into remanufacturing models.

\begin{table}[ht!]
\centering
\caption{Summary of Related Literature}\label{tab:literature}
 \resizebox{\textwidth}{!}{
\begin{tabular}{lccccc}
    \toprule
    \multirow{2}{*}{} & \multicolumn{2}{c}{OEM's Business Model} & \multirow{2}{*}{Consumer perception} & \multirow{2}{*}{Market size} & \multirow{2}{*}{Contract Form}\\
    \cline{2-3}
    & In-house & Cooperation & & & \\
    \midrule
     \cite{agrawal2010essays} & \checkmark &$\times$ & \checkmark  & Constant & $\times$\\
     
     \cite{fang2020third} & \checkmark &$\times$ &  \checkmark & Constant& $\times$\\         \cite{huang2024remanufacturing} & \checkmark & \checkmark & $\times$ & Constant & One-part\\
     
     \cite{li2024should} & $\times$ & \checkmark & \checkmark  & Constant & One-part\\
     
     This work & \checkmark & \checkmark & \checkmark  & Stochastic& Two-part Tariff\\
    \bottomrule
\end{tabular}
}
\end{table}

Table~\ref{tab:literature} compares related studies and highlights our contributions. This work contributes to the literature by jointly considering consumer perceptions, stochastic market size, and a two-part authorization contract. 

\section{Model formulation} \label{sec:model framework}
Our primary objective is to determine the most profitable remanufacturing model for OEM under different consumer perception conditions. To address this, we consider three alternative models: no remanufacturing (Model N), OEM in-house remanufacturing (Model O), and TPR-authorized remanufacturing (Model T). 

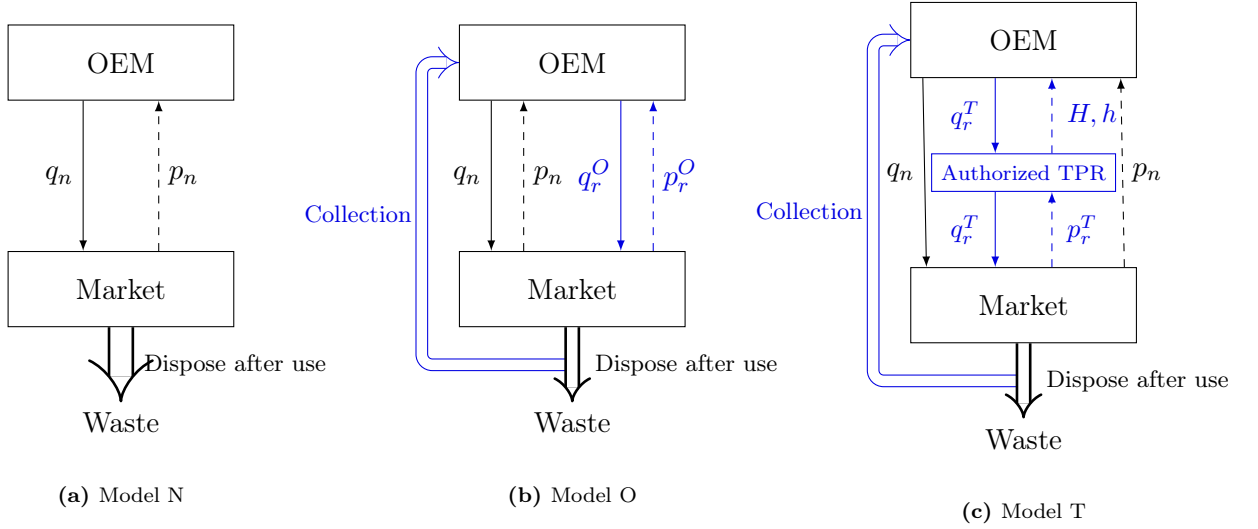
\begin{figure}[ht!]
\definecolor{reman_color}{RGB}{0,0,225} 

\begin{tikzpicture}[
        node distance=2cm, 
        auto,
        box/.style={draw, rectangle, minimum width=3cm, minimum height=1cm},
        thick arrow/.style={-implies, double, double distance=3pt, line width=1pt},
        small label/.style={font=\footnotesize}
    ]
% Subfigure 1
\begin{scope}[yshift=-0.3cm]
            \node[box] (oem1) {OEM};
            \node[box] (market1) [below=of oem1] {Market};
            \node[below=1cm of market1] (waste1) {Waste};
            \draw[-latex] ($(oem1.south)!1/3!(oem1.south west)$) -- node[left] {$q_n$} ($(market1.north)!1/3!(market1.north west)$);
            \draw[-latex, dashed] ($(market1.north)!1/3!(market1.north east)$) -- node[right] {$p_n$} ($(oem1.south)!1/3!(oem1.south east)$);
            \draw[-implies, double, double distance=8pt, line width=1pt, minimum width=3cm] (market1.south) -- node[right, small label] {Dispose after use} (waste1.north);
            \node [below=0.1cm of waste1] {\parbox{0.3\linewidth}{\subcaption{Model N}}};
    \end{scope}

        % Subfigure 2
\begin{scope}[xshift=6cm, yshift=-0.3cm]
            \node[box] (oem2) {OEM};
            \node[box] (market2) [below=of oem2] {Market};
            \node[below=1cm of market2] (waste2) {Waste};
            
            \draw[-latex] ($(oem2.south)!5/7!(oem2.south west)$) -- node[left] {$q_n$} ($(market2.north)!5/7!(market2.north west)$);
            \draw[-latex, dashed] ($(market2.north)!3/7!(market2.north west)$) -- node[right] {$p_n$} ($(oem2.south)!3/7!(oem2.south west)$);
            \draw[-latex, reman_color] ($(oem2.south)!3/7!(oem2.south east)$) -- node[left] {$q_r^O$} ($(market2.north)!3/7!(market2.north east)$);
            \draw[-latex, dashed, reman_color] ($(market2.north)!5/7!(market2.north east)$) -- node[right] {$p_r^O$} ($(oem2.south)!5/7!(oem2.south east)$);
            \draw[-implies, reman_color, rounded corners, double, double distance=4pt] 
    ($(waste2.north)!0.5!(market2.south)$) -- ++(-2,0) |- 
    node[left, reman_color, pos=0.25, small label] {Collection} (oem2.west);
            \draw[-implies, double, double distance=4pt, line width=1pt, minimum width=3cm] (market2.south) -- node[right, small label] {Dispose after use} (waste2.north);
            \node [below= 0.1cm of waste2] {\parbox{0.3\linewidth}{\subcaption{Model O}}};
    \end{scope}       
    
        % Subfigure 3
    \begin{scope}[xshift=12cm]
            \node[box] (oem3) {OEM};
            \node[box, reman_color, minimum width=1cm, minimum height=0.5cm, font = \scriptsize] (tpr3) [below=1cm of oem3] {Authorized TPR};
            \node[box] (market3) [below=1cm of tpr3] {Market};
            \node[below=1cm of market3] (waste3) {Waste};
            \draw[-latex] ($(oem3.south)!8/9!(oem3.south west)$) -- node[left] {$q_n$} ($(market3.north)!6/7!(market3.north west)$);
            \draw[-latex, dashed] ($(market3.north)!8/9!(market3.north east)$) -- node[right] {$p_n$} ($(oem3.south)!6/7!(oem3.south east)$);
            \draw[-latex, reman_color] ($(oem3.south)!1/4!(oem3.south west)$) -- node[left=2pt] {\small $q_{r}^T$} ($($(oem3.south)!1/4!(oem3.south west)$)+(0,-1)$);
            \draw[-latex, dashed, reman_color] ($($(oem3.south)!1/4!(oem3.south east)$)+(0,-1)$) -- node[right=2pt] {\small $H, h$} ($(oem3.south)!1/4!(oem3.south east)$);
            \draw[-latex, reman_color] ($($(market3.north)!1/4!(market3.north west)$)+(0,1)$) -- node[left=2pt]  {\small $q_r^T$} ($(market3.north)!1/4!(market3.north west)$);
            \draw[-latex, dashed, reman_color] ($(market3.north)!1/4!(market3.north east)$) -- node[right=2pt] {\small $p_r^T$} ($($(market3.north)!1/4!(market3.north east)$)+(0,1)$);
            \draw[-implies, reman_color, rounded corners, double, double distance=4pt] 
    ($(waste3.north)!0.5!(market3.south)$) -- ++(-2,0) |- 
    node[left, reman_color, pos=0.25, small label] {Collection} (oem3.west);
            \draw[-implies, double, double distance=4pt, line width=1pt, minimum width=3cm] (market3.south) -- node[right, small label] {Dispose after use} (waste3.north);
        
    \node [below = 0.1cm of waste3] {\parbox{0.3\linewidth}{\subcaption{Model T}}};
    \end{scope}
    \end{tikzpicture}
\caption{Model Structure} \label{fig:model structures}
\end{figure}

Let $p_n$ and $p_r$ denote the prices of new and remanufactured products, respectively, and $q_n$ and $q_r$ denote their corresponding quantities. Superscripts indicate the decision-maker, with \textit{O} referring to the OEM and \textit{T} to the TPR. For notational clarity, these superscripts will be omitted when the context of the remanufacturing model is unambiguous. Fig.~\ref{fig:model structures} presents an overview of the decision-making frameworks considered. 

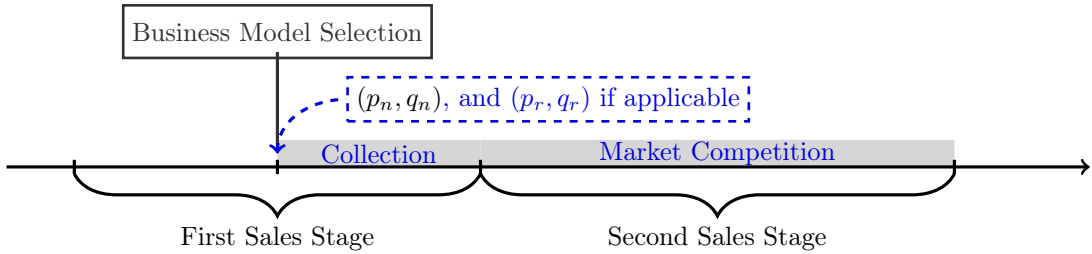
\begin{figure}[ht!]
\centering
\definecolor{ColorOne}{RGB}{50, 50, 50}
\definecolor{ColorTwo}{RGB}{150, 150, 150} 
\definecolor{ColorThree}{RGB}{0,0,225} 

\tikzstyle{descript} = [text = black, align=center, minimum height=1cm, align=center, outer sep=0pt,font = \footnotesize]
\tikzstyle{activity} =[align=center,outer sep=0pt]

\begin{tikzpicture}[very thick, black,scale=0.9]
\small

%% Coordinates
\coordinate (O) at (-1,0); % 
\coordinate (S) at (0,0); % start of sales
\coordinate (1E) at (6,0); % end of 1st sales
\coordinate (2E) at (13,0); % end of 2nd sales
\coordinate (F) at (15,0); %End
\coordinate (E1) at (3,0); %business decisions

% First and Second Sales Stages
\draw [decorate,decoration={brace,amplitude=16pt,mirror}]($(S)+(0,-0.1)$) -- ($(1E)+(0,-0.1)$) node [black,midway,below=16pt] {First Sales Stage};
\draw [decorate,decoration={brace,amplitude=16pt,mirror}]($(1E)+(0,-0.1)$) -- ($(2E)+(0,-0.1)$) node [black,midway,below=16pt] {Second Sales Stage};

%% if remanufacture
% \path[dashed, <-,color=ColorOne] ($(E1)+(0,0.2)$) edge [out=90, in=-100] ($(E1)+(1.5,1)$);
% \draw ($(E1)+(0,0.2)$) node[below=0pt,align=left,ColorThree] {{\color{black}$(p_n,q_n)$} \\ $(p_r,q_r)$, if applicable};

% collecting process
\fill[color=ColorOne!20] rectangle ($(E1)$) -- ($(1E)$) -- ($(1E)+(0,0.4)$) -- ($(E1)+(0,0.4)$);
\draw ($(E1)+(1.5,0.2)$) node[activity,ColorThree] {Collection};

% Reman. sales
\fill[color=ColorOne!20] rectangle ($(1E)$) -- ($(2E)$) -- ($(2E)+(0,0.4)$) -- ($(1E)+(0,0.4)$);
\draw ($(2E)+(-3.5,0.2)$) node[activity,ColorThree] {Market Competition};
% decide price and quantity
% \node[draw=ColorOne,text=ColorOne, minimum height=0.6cm, minimum width= 1.2cm, line width=1pt, align=left](D3) at ($(E1)+(-1,1)$) {$(p_n,q_n)$};
% \path[->,color=ColorOne] ($(D3)+(0,-0.3)$) edge [out=-50, in=180-50] ($(E1)+(0,0.25)$);

%% Select Business Model
\node[draw=ColorOne,text=ColorOne, minimum height=0.7cm, minimum width= 1cm, line width=1pt, align=center](D1) at ($(E1)+(0,2)$) {Business Model Selection};
\path[<-,color=ColorOne] ($(E1)+(0,0.2)$) edge [out=90, in=-90]  ($(D1)+(0,-0.3)$);

% decide price and quantity
\node[dashed,draw=ColorThree,text=ColorThree, minimum height=0.3cm, minimum width= 1cm, line width=1pt, align=left](D2) at ($(E1)+(4,1)$) {{\color{black}$(p_n,q_n)$}, and $(p_r,q_r)$ if applicable};
\path[dashed, <-,color=ColorThree] ($(E1)+(0,0.2)$) edge [out=90, in=180] (D2.west);

%% Arrow
\draw[->] (O) -- (F);
%% Ticks
\foreach \x in {0,3,6,13}
\draw(\x cm,3pt) -- (\x cm,-3pt);
%% Labels
% \foreach \i \j in {0/T=0,6/T=$t_1$,13/T=$t_1+t_2$}{
% 	\draw (\i,0) node[below=3pt] {\j} ;
% }

\end{tikzpicture}
\caption{Timeline} \label{fig:timeline}
\end{figure}

The sequence of events is as follows. In the first sales stage, the OEM decides whether to engage in remanufacturing. If the OEM opts against remanufacturing, only new products are offered in the market. If the OEM opts to remanufacture, the OEM begins to collect the used products as remanufacturing materials. This practice is common in industry, as seen in trade-in rebate programs by Philips and Apple, or Dell's take-back program, which allows customers to exchange used products for credits or discounts on future purchases. The OEM also determines who will perform remanufacturing, and relevant parties decide on pricing and production quantities for their respective products. In the second sales stage, if remanufacturing is performed, remanufactured products compete with new products in the market. 
Fig.~\ref{fig:timeline} illustrates the timeline of the entire process. 

It is worth noting that our analysis concentrates on the market competition between new and remanufactured products at the second stage, rather than the dynamic interactions across both stages. We aim to maximize the OEM's expected profit at the second stage; hence, we have a single-period problem.
The optimization problem is formulated as follows. 
\begin{equation*}
\max_{\mathbbm{1}_N,\mathbbm{1}_O,\mathbbm{1}_T} \Pi_{OEM} = \mathbbm{1}_N\Pi^{N}(p_n, q_n) + \mathbbm{1}_O\Pi^{O}(p_n, q_n, p_r^O, q_r^O) + \mathbbm{1}_T\Pi^{T}(p_n, q_n, p_r^T, q_r^T).
\label{eq:general profit}
\end{equation*}

If there is no remanufacturing (Model N), the OEM's expected profit equals the revenue from new products minus production costs for new units: $\Pi^{N}(p_n, q_n) = \mathbb{E}\left[ p_n S_n - c q_n \right]$, which follows the classical newsvendor problem. For OEM in-house remanufacturing (Model O), expected profit is revenues from new and remanufactured products, minus production, remanufacturing, and collection costs: $\Pi^{O}(p_n, q_n,p_r, q_r) = \mathbb{E}\left[ p_n S_n - c q_n + p_r^O S_r - c_r q_r^O - c_{coll} q_r^O \right]$. For TPR-authorized remanufacturing (Model T), the OEM earns revenue from new products and authorization fees from the TPR, while incurring production and collection costs: $\Pi^{T}(p_n, q_n, p_r^T, q_r^T) = \mathbb{E}\left[ p_n  S_n - c  q_n + H + h q_r^T - c_{coll} q_r^T\right]$. Here, $S_n$ and $S_r$ denote sales of new and remanufactured products. $\mathbb{E}[\cdot]$ and $\mathbbm{1}_{[\cdot]}$ denote the expectation and indicator operators, respectively. Table~\ref{tab:notations} summarizes our notation.
\begin{table}[ht!]\small
\centering
\caption{Notation}\label{tab:notations}
\begin{tabular}{ll}
        \toprule
        \textbf{Exogenous variables} & \textbf{Definition} \\
        \midrule
        $V/V_n/V_r$ & Base/ New/ Remanufactured product value\\
        $\theta \sim U[0,1]$ & Consumer preference, uniformly distributed in $[0,1]$\\
        $\delta \in [0,1]$ & Depreciation factor for products sold in second stage \\
        $\alpha \in [0,1]$ & WTP discount factor for remanufactured products, $V_r = \alpha  V_n$ \\
        $\beta\in [-1,1]$ & Perception factor of  assimilation/ contrast effect\\
        $c$ & New production cost\\
        $c_r$ & Remanufacturing cost \\
        $c_{\text{coll}}$ & Collection cost \\
        $N \sim \text{Poisson}(\lambda)$ & Stochastic market size, with fixed expected size $\lambda$ \\
        \midrule
      \textbf{Endogenous variables} & \textbf{Definition} \\
        \midrule
        $U_{n}/U_r$ & Consumer utility from new/remanufactured product \\
        $D_n/D_r$ & Demand for new/remanufactured product\\
        $S_n/S_r$ & Sales of new/remanufactured product\\
        $H/h$ & One-time/ unit authorization fee\\
        $\Pi^N/\Pi^O/\Pi^T$ & OEM's expected profit under Model N/O/T \\
        \midrule
        \textbf{Decision Variables} & \textbf{Definition} \\
        \midrule
        $p_n/p_r^O/p_r^T\in \mathbb{R^+}$ & Price of new/ OEM-remanufactured/ TPR-remanufactured product\\
        $q_n/q_r^O/q_r^T\in \mathbb{N}$ & Quantity of new/ OEM-remanufactured/ TPR-remanufactured products \\
        \bottomrule
    \end{tabular}
\end{table}

Next, we detail the consumer utility with perception effects. We consider a market in which the total number of potential customers follows a Poisson distribution with a fixed rate, that is, $N\sim \text{Poisson}(\lambda)$. Each customer is characterized by a heterogeneous preference $\theta \sim U[0,1]$. Their willingness-to-pay (WTP) for a product depends on both $\theta$ and the value of the product. New products have a base value $V$ in the first stage. In the second stage, new products depreciate to $V_n=\delta V$ with $0 \leq \delta \leq 1$. Remanufactured products are perceived as less valuable than new ones; this is captured by a discount factor $\alpha\in[0,1]$, such that their value is given by $V_r=\alpha V_n$. 
As discussed in Section \ref{sec:Literature review}, the perceived value of new products is influenced by the identity of the remanufacturer through assimilation and contrast effects \citep{agrawal2015remanufacturing}. To capture these effects, we adjust the perceived value of new products when both products coexist in the market to
$g(V_n, V_r)= V_n + \beta (V_n - V_r)$, where $-1\leq\beta\leq 1$ quantifies the magnitude of assimilation and contrast effects. In particular, $\beta^-<0$ reflects the assimilation effect under Model O, and $\beta^+>0$ reflects the contrast effect under Model T. Thus, the consumer's WTP for a new product is adjusted to $\theta g(V_n, V_r)$, while the WTP for a remanufactured product is $\theta V_r$. Consumer utility is calculated as the difference between the WTP and the price of the product. When new and remanufactured products are priced at $p_n$ and $p_r$, respectively, the utility of a new product is $U_n=\theta g(V_n, V_r) -p_n$, and for a remanufactured product is $U_r=\theta V_r-p_r$. Customers buy the product that offers the highest nonnegative utility; in particular, if $U_n=U_r$, consumers default to choosing the new product. 

\section{Price and quantity decisions under remanufacturing business models} \label{sec:remanufacturing modes}
In this section, we develop the market demand and expected profit functions for each of the three business models. For every model, we derive the OEM's profit-maximizing prices and quantities for both new and remanufactured products.

\subsection{Model N: no remanufacturing}

To establish a baseline, we first analyze the OEM's profit without remanufacturing. 
In this model, the OEM sets the price and production quantity to maximize the expected profit $\Pi^{N}(p_n, q_n) = \mathbb{E}\left[ p_n S_n - c q_n \right]$, where $c$ denotes the unit production cost and $S_n = \min\{D_n,q_n\}$ denotes the sales volume. This corresponds to the classical price-setting newsvendor problem. 

Before defining the optimization problem, we first characterize the demand $D_n$. A customer with preference $\theta$ buys the new product if the utility is nonnegative, i.e, $U_n(\theta)=\theta V_n-p_n\geq 0$. Thus, the demand consists of all consumers with $\theta \geq \frac{p_n}{V_n}$. Given a Poisson-distributed market size, i.e., $N\sim\text{Poisson}(\lambda)$ and uniform consumer preferences, i.e., $\theta\sim U[0,1]$, market demand $D_n$ follows a compound Poisson distribution with a demand rate $\Lambda(p_n):= (1-\frac{p_n}{V_n}) \lambda$. The cumulative distribution function (CDF) for the demand $D_n$ is $F(k,\Lambda(p_n))=\sum_{i=0}^k \frac{1}{i!}e^{-\Lambda(p_n)} \Lambda^i(p_n)$. All detailed proofs can be found in the Appendix.

Substituting the demand rate into the expected profit function, we formulate the optimization problem as
\begin{equation*}
\max_{p_n, q_n} \ \Pi^N(q_n,p_n) = p_n e^{-\Lambda(p_n)} \left[\sum_{k=1}^{q_n}\frac{\Lambda^k(p_n)}{(k-1)!} + \sum_{k=q_n+1}^{\infty} q_n \frac{\Lambda^k(p_n)}{k!} \right] -c q_n.
\label{eq:model n optimization problem}
\end{equation*}

Following classical price-setting newsvendor theory, we first solve the optimal quantity at a given price, then substitute this value into the expected profit function to optimize over price. For any $p_n\in (c,V_n)$, the OEM's optimal quantity is characterized by the critical fractile $q_n^*(p_n)=F^{-1}\Big(1-\frac{c}{p_n}\Big)$, which is consistent with the literature \citep{arrow1951optimal,littlewood1972forecasting}. Here, $F^{-1}(y):=\min{\{k \in \mathbb{N}:F(k,\Lambda(p_n))\geq y}\}$ denotes the inverse CDF. Restricting $p_n\in (c,V_n)$ excludes trivial cases where the OEM might theoretically produce zero or infinite quantities, as in the standard literature \citep{gallego1994optimal,petruzzi1999pricing}. Substituting $q_n^*(p_n)$ into the expected profit function generates the reduced-form pricing problem:
\begin{equation}
\max_{p_n} \ \Pi^N(p_n,q_n^*)= p_n \Lambda(p_n) F(q_n^* -1,\Lambda(p_n)).
\label{eq: model n expected profit with optimal q}
\end{equation}

% approximation 

Solving this pricing problem analytically is challenging because price $p_n$ is continuous while production quantity $q_n$ is an integer. This creates a ``saw-tooth'' pattern in the objective function $\Pi^N(p_n,q_n^*)$, arising from the discontinuous jumps in the critical fractile solution $q_n^*(p_n)=F^{-1}(1-\frac{c}{p_n})$ as $p_n$ varies continuously. 
% Although \citet{schulte2020price} analyze a similar Poisson demand problem, their approach relies on computationally intensive enumeration. 
Therefore, we use a numerical approach to solve the problem.

\subsection{Model O: OEM in-house remanufacturing} 

In Model O, the OEM produces both new and remanufactured products. The OEM sets the prices and quantities for new products $(p_n,q_n)$ and remanufactured products $(p_r,q_r)$. To remanufacture, the OEM first collects the used products at a unit collection cost of $c_{coll}$. All collected units are subsequently remanufactured, so the total collection quantity is $q_r$. The unit cost to produce a new product is $c$ and the unit cost of remanufacturing (excluding collection) is $c_r$. The unit cost form is widely adopted in the remanufacturing literature, such as \citet{atasu2008remanufacturing}. The OEM's expected profit function is $\Pi^{\text{O}}(p_{n},q_{n},p_r,q_r) = \mathbb{E}\big[ p_{n} S_n + p_r   S_r - c q_{n}- (c_r+c_{coll}) q_r \big]$, where $S_n:=\min \{D_n,q_n\}$ denotes the sales of new products and $S_r:=\min\{D_r,q_r\}$ denotes the sales of remanufactured products. The goal is to determine the prices and quantities for both products so that the expected profit is maximized.

To determine the optimal decisions of the OEM, we first derive the demands for both products. According to \citet{agrawal2015remanufacturing}, when both new and remanufactured products are offered by the same manufacturer, i.e., the OEM, consumers tend to perceive their values as similar, which reduces their perceived value of new products. This effect, the so-called \textit{assimilation effect}, is captured by $\beta^-\in [-1,0]$. Each customer evaluates utility of the new product $U_n=\theta V_n + \theta \beta^- (V_n - V_r) -p_n$ and of the remanufactured product $U_r=\theta V_r-p_r$, where $V_n = \delta V$, $V_r = \alpha \delta V$. A customer buys a new product if the utility of the new product is nonnegative and is at least as high as the utility of the remanufactured product; that is, if $U_n\geq 0$ and $U_n\geq U_r$. The customer buys a remanufactured product if the utility of the remanufactured product is nonnegative and is greater than the utility of the new product; that is, if $U_r\geq0$ and $U_r>U_n$. If neither condition is met, the customer does not buy. Thus, the demand for each product is determined by the number of customers whose preferences satisfy the respective utility condition. 

Given a Poisson-distributed market size ($N\sim\text{Poisson}(\lambda)$) and uniform customer preferences ($\theta\sim U[0,1]$), the demands for both new products $D_n$ and remanufactured products $D_r$ follow compound Poisson distributions with CDFs $F_n(k,\Lambda_n^O(p_n,p_r))$ and $F_r(k,\Lambda_r^O(p_n,p_r))$, respectively. Their respective demand rates are as follows: 
\begin{align}
    &\Lambda_{n}^O(p_n,p_r):= \left(1-\min\left\{1,\max\left\{\frac{p_{n}-p_r}{(1+\beta^-)(1-\alpha)\delta V }, \frac{p_{n}}{(1+\beta^- -\alpha \beta^-)\delta V }\right\}\right\}\right) \lambda, \label{eq: Model O lambda n}\\
    &\Lambda_{r}^O(p_n,p_r):= \left( \max\left\{ 0, \min\left\{1,\frac{p_{n}-p_r}{ (1+\beta^-)(1-\alpha)\delta V}\right\} - \min\left\{1,\frac{p_{r}}{\alpha \delta V}\right\} \right\}\right) \lambda. \label{eq: Model O lambda r}
\end{align} 

\begin{figure}[ht!]
\centering
\includegraphics[width = 0.8\textwidth]{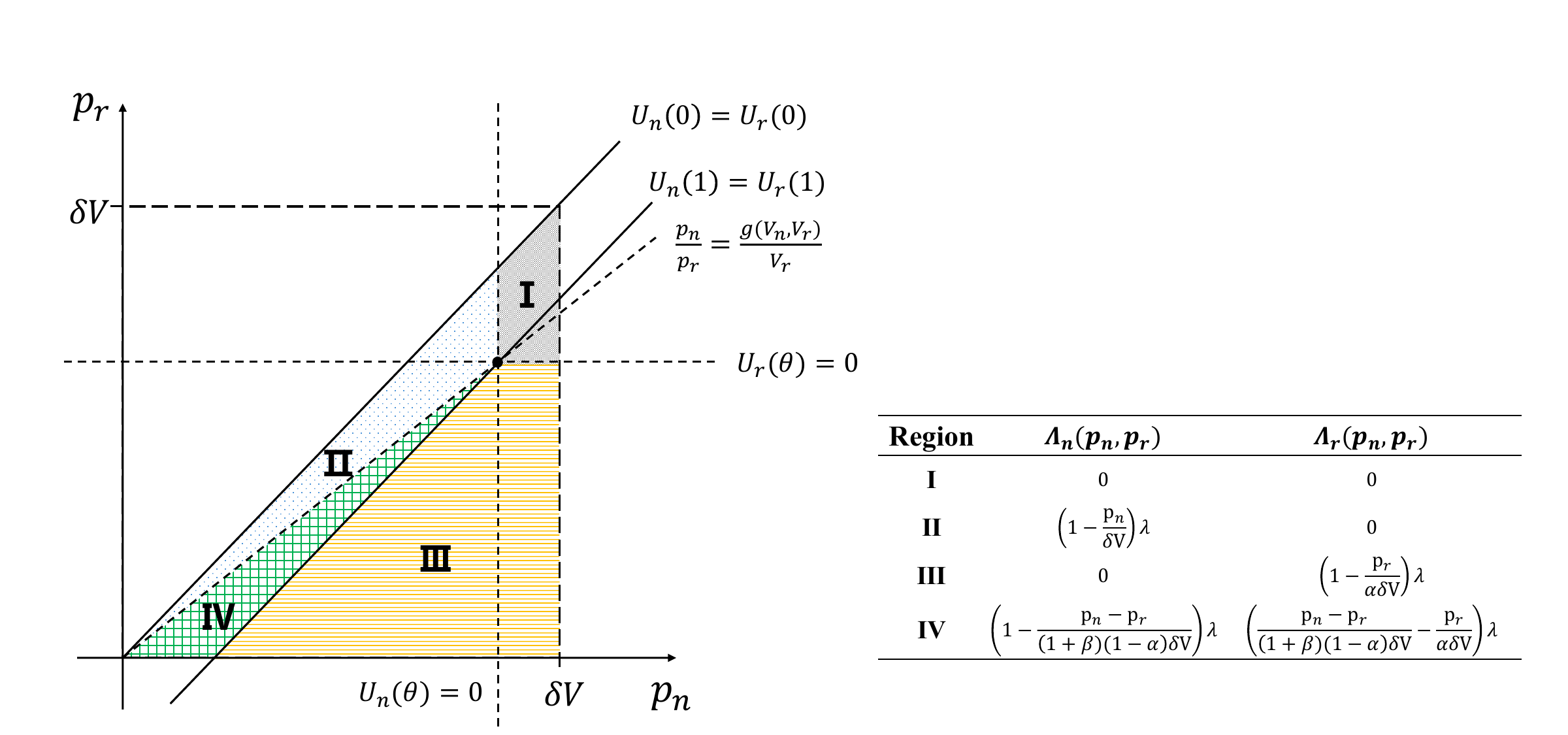}
\caption{Region chart of $\Lambda_n$ and $\Lambda_r$ with $\beta^- \in [-1,0]$} \label{fig: Model O demand rates cases}
\end{figure}

Fig.~\ref{fig: Model O demand rates cases} illustrates four possible market outcomes and the joint value of $\Lambda_n^O(p_n,p_r)$ and $\Lambda_r^O(p_n,p_r)$ as functions of $p_n$ and $p_r$. For simplicity, in this section we abbreviate $\Lambda_n^O(p_n,p_r),\Lambda_r^O(p_n,p_r)$ as $\Lambda_n,\Lambda_r$, respectively. 
Contrary to prior studies that assume the natural coexistence of new and remanufactured products on the market \citep{atasu2008remanufacturing,xiong2013don,zou2016third}, we show a wider range of market structures: not only coexistence, but also scenarios where remanufacturing is infeasible (Region \uppercase\expandafter{\romannumeral1} and \uppercase\expandafter{\romannumeral2}) and where remanufactured products dominate the market (Region \uppercase\expandafter{\romannumeral3}). 

Substituting the demand rates into the expected profit function, we formulate the optimization problem as
\begin{align*}
\max_{p_{n},q_{n},p_r,q_r} \ \Pi^{\text{O}}(p_{n},q_{n},p_r,q_r) = \ 
    & p_{n}e^{-\Lambda_{n}} \left[\sum_{k=1}^{q_{n}} \frac{1}{(k-1)!} \Lambda_{n}^k+ \sum_{k=q_{n}+1}^{\infty} q_{n} \frac{1}{k!} \Lambda_{n}^k\right] \notag \\
    &+ p_r e^{-\Lambda_{r}} \left[\sum_{k=1}^{q_{r}} \frac{1}{(k-1)!} \Lambda_{r}^k + \sum_{k=q_{r}+1}^{\infty} q_{r} \frac{1}{k!} \Lambda_{r}^k\right] \notag \\
    & -c q_{n}-(c_r+c_{coll}) q_r.
\label{eq: Model O profit}
\end{align*}

By solving the first- and second-order conditions, we obtain the optimal quantities for new products as $q_{n}^{*}(p_n,p_r)=F_n^{-1}\left(1-\frac{c}{p_{n}}\right)$ and for remanufactured products as $q_{r}^{*}(p_n,p_r)=F_r^{-1}\left(1-\frac{c_r+c_{coll}}{p_{r}}\right)$, for any given prices $p_n\in(c,V_n)$ and $p_r\in(c_r+c_{coll},V_r)$. 
Here, $F_n^{-1}$ and $F_r^{-1}$ are the inverse CDFs of demands for new and remanufactured products. For simplicity, in this section we abbreviate optimal quantities $q_{n}^{*}(p_n,p_r), q_{r}^{*}(p_n,p_r)$ as $q_n^*,q_r^*$, respectively. Substituting the optimal quantities into the expected profit function, we reformulate the optimization problem
\begin{equation}
\max_{p_n,p_r}\ \Pi^{\text{O}}(p_{n},p_r) = p_{n}\Lambda_{n} F_n\left(q^*_{n}-1,\Lambda_n\right)+ p_r \Lambda_r F_r\left(q_r^*-1,\Lambda_r\right).
\label{eq: Model O profit with optimal q}
\end{equation}

% approximation
As in Model N, the profit function in Model O exhibits a ``3D saw-tooth'' pattern due to the combination of continuous prices and discrete quantities. We use a numerical method to solve this optimization problem in the next section.

% In the supplementary material, we analyze the impact of $\alpha$ and $\beta^-$ on $\widetilde{\Pi}^{O*}$ and examine the extreme cases with $\alpha=0,1,\beta^-=-1,0$.

\subsection{Model T: TPR-authorized remanufacturing}

In Model T, the OEM authorizes a TPR to remanufacture products through a two-part tariff contract. The OEM is responsible for the production and sales of new products and the collection of used items, incurring a unit production cost $c$ and a unit collection cost $c_{coll}$. The OEM decides the price and quantity of new products ($p_n,q_n$). The authorized TPR is responsible for the remanufacturing and sales of remanufactured products, incurring unit remanufacturing cost $c_r$. The TPR decides the price and quantity of remanufactured products ($p_r,q_r$). The authorization contract requires TPR to pay a fixed one-time fee $H$ for technology permission and a unit fee $h$ for each remanufactured item. Restricting $h\in(0,c)$ prevents the TPR from having a financial incentive to artificially create new products only for the purpose of remanufacturing them later. This ensures that TPR engages in remanufacturing primarily with legitimately used products. In our analysis, the contract parameters $(H,h)$ are treated as exogenously given, but in Section~\ref{sec:contract and stochasticity}, we perform a sensitivity analysis on the value of $(H,h)$. 

The interaction between the OEM and the TPR is formulated as a Stackelberg game, with the OEM as the leader and the TPR as the follower. The OEM first determines $p_n$ and $q_n$. Observing these decisions, the TPR chooses the optimal $p_r$ and $q_r$ to maximize its expected profit $\Pi_{tpr}^T$. The cooperation of authorization proceeds only when the TPR's maximized profit is nonnegative.
The OEM aims to maximize its expected profit $\Pi_{oem}^T$ by selecting the optimal $(p_n, q_n)$ through anticipating the optimal response of the TPR. Given the setting of Stackelberg game, we formulate the optimization problem of Model T as
\begin{align*}
\max_{p_n,q_n}\quad &\Pi_{oem}^T (p_n,q_n,p_r,q_r):= \mathbb{E}\big[p_{n} S_n -c q_{n} + H + (h -c_{coll}) q_r\big] \\
\text{s.t.}\quad & \Pi_{tpr}^T(p_n,q_n,p_r^*,q_r^*) := \mathbb{E}\big[p_{r} S_r- H - (h + c_r)  q_r\big] \geq0, \\
& (p_r^*,q_r^*)\in \argmax_{p_r,q_r}\ \Pi_{tpr}^T(p_n,q_n,p_r,q_r).
\end{align*}

Analogously to Model O, we start by deriving the demands for both products. Each customer compares the utility of the new product, i.e., $U_n(\theta)=\theta V_n+\theta\beta^+(V_n-V_r)-p_n$, and of the remanufactured product, i.e., $U_r(\theta)=\theta V_r-p_r$, where $V_n = \delta V$ and $V_r = \alpha \delta V$. The parameter $\beta^+\in[0,1]$ captures the \textit{contrast effect}: when new and remanufactured products are offered by different parties, i.e., the OEM and TPR, respectively, customers tend to perceive a strong distinction between the two products and enhance the perceived value of new products \citep{agrawal2015remanufacturing}. Customers buy the product that leads to highest nonnegative utility.

The demand rates $\Lambda_n^T(p_n,p_r)$ and $\Lambda_r^T(p_n,p_r)$ take the same form as Eqs.~\eqref{eq: Model O lambda n}-\eqref{eq: Model O lambda r}, but with $\beta^+$ instead of $\beta^-$.
% Denote $F_n(k,\Lambda_{n}^T(p_n,p_r))$ and $F_r(k,\Lambda_{r}^T(p_n,p_r))$ as CDFs for demands for new $D_n$ and remanufactured $D_r$ products, respectively. 
For simplicity, in this section we abbreviate $\Lambda_{n}^T(p_n,p_r)$ and $\Lambda_{r}^T(p_n,p_r)$ to $\Lambda_{n}$ and $\Lambda_{r}$, respectively.
\begin{figure}[ht!]
\centering
\includegraphics[width = 0.8\textwidth]{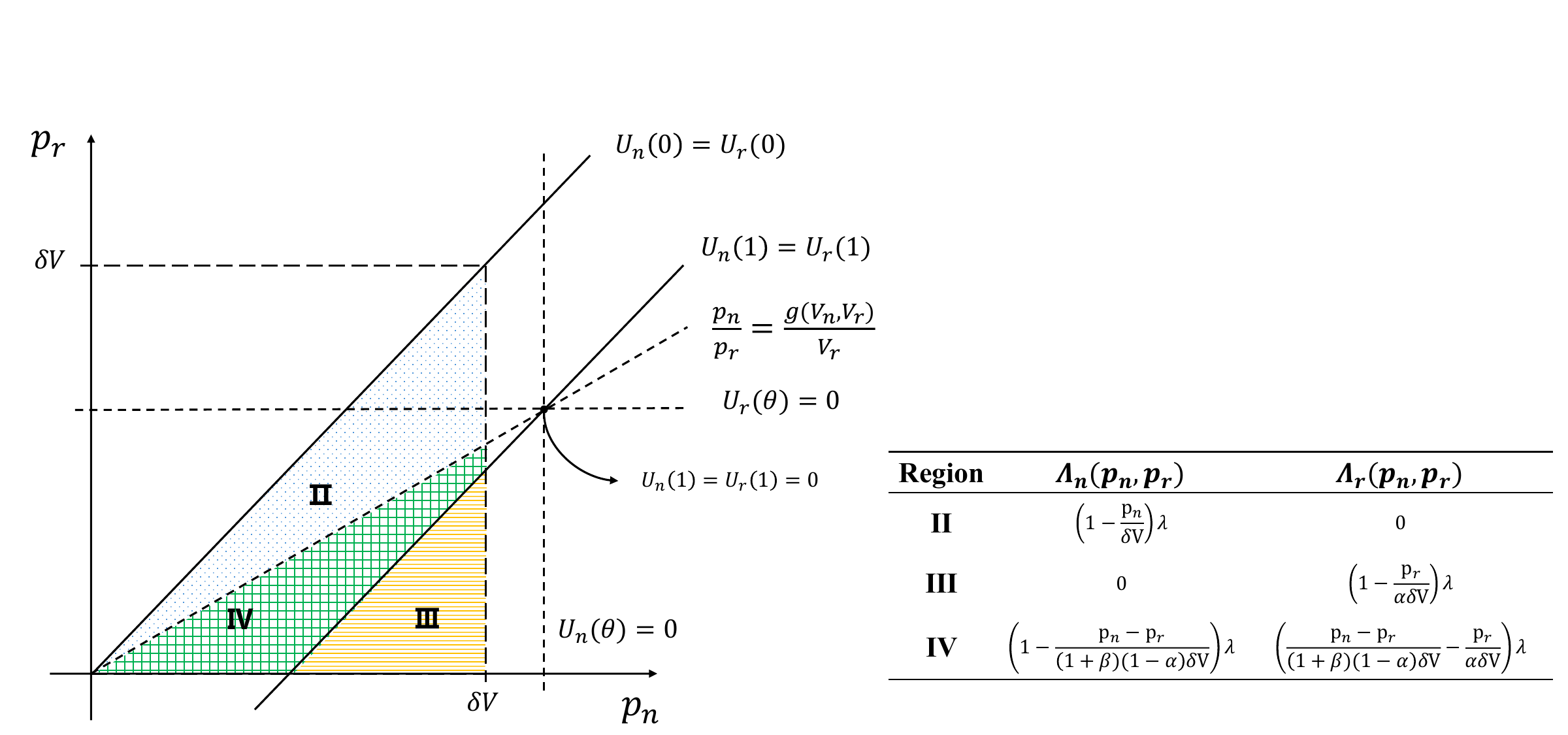}
\caption{Region chart of $\Lambda_n$ and $\Lambda_r$ with $\beta^+ \in [0,1]$} \label{fig: Model T regions}
\end{figure}

Fig~\ref{fig: Model T regions} illustrates three possible market outcomes in Model T, along with the joint values of $\Lambda_n$ and $\Lambda_r$ for each outcome. Contrary to Model O, Model T has only three distinct regions: the only new market (Region \uppercase\expandafter{\romannumeral2}), the only remanufactured market (Region \uppercase\expandafter{\romannumeral3}), and the coexistence market (Region \uppercase\expandafter{\romannumeral4}). The insight is that consumers with high preferences always consider new products due to the contrast effect, i.e., $U_n(\theta)\geq 0$ for all $\theta$. 

Applying backward induction, the OEM anticipates that for any price $p_n$ of new products, the TPR selects its optimal price and quantity for the remanufactured product, that is, $p^*_r(p_n)$ and $q_{r}^{*}(p_r^*)=F_r^{-1}(1-\frac{c_r+h}{p_{r}^*})$. The OEM then chooses $p_n$ and $q^*_n(p_n)=F_n^{-1}(1-\frac{c}{p_{n}})$ to maximize its profit. Upon now, there are two main challenges to this sequential optimization. First, because the discrete quantity $q^*_r$ is isolated rather than being embedded within a CDF in $\Pi_{oem}^T(p_n,p_r^*)$, it is challenging to obtain a closed-form solution for the OEM optimization problem. Second, the contract parameters $(H,h)$ further complicates the mathematics and prevents an explicit solution. As a result, we solve the model by numerical method in the next section.

\section{Results and analysis} \label{sec:selection}
In this section, we identify the optimal business model across a range of consumer perception conditions, using numerical approach. We translate these results into a hierarchical decision roadmap that guides the selection of remanufacturing business models. Finally, we assess the resulting market dynamics and environmental impact.

\subsection{Selection of remanufacturing business models}

We apply numerical methods to systematically compare the profitability of all three remanufacturing business models under all range of consumer perception conditions, i.e., $\alpha\in[0,1]$ and $|\beta|\in[0,1]$. Following parameter setting from \citet{ferrer2006managing} and related literature, we set the following parameter values in our numerical studies: expected market size $\lambda =1000$, base value of the product $V=1000$, depreciation factor $\delta=0.8$, unit production cost for the new product $c=200$, unit collection cost for used products $c_{coll}=40$ and unit remanufacturing cost $c_r=80$. The total unit cost of remanufacturing a product $c_r+c_{coll}=120$, which is consistent with empirical and industry evidence that remanufacturing typically reduces costs by 40\% - 65\% \citep{ginsburg2001once,savaskan2004closed,du2012integrated}. The authorization contract specifies a fixed one-time fee $H=10,000$ and a unit fee $h=100$, and these terms $(H,h)$ will be discussed further in the next section.

For each model, we use a grid search to find prices and quantities that maximize the OEM's expected profit, using a price increment of $0.01$ within their respective support set. For each pair of prices $(p_n,p_r)$, the optimal quantities are set as the lower integer of the critical fractile, then these values are used to compute the expected profit. The OEM then selects the business model with the highest maximum expected profit. We repeat this procedure over $\alpha\in[0,1]$ and $|\beta|\in[0,1]$ with an increase of $0.01$. For convenience, we assume the same degree of assimilation and contrast effects, i.e., $\vert \beta^+ \vert=\vert \beta^- \vert$.

\begin{figure}[h!]
\centering
\includegraphics[width=0.7\textwidth]{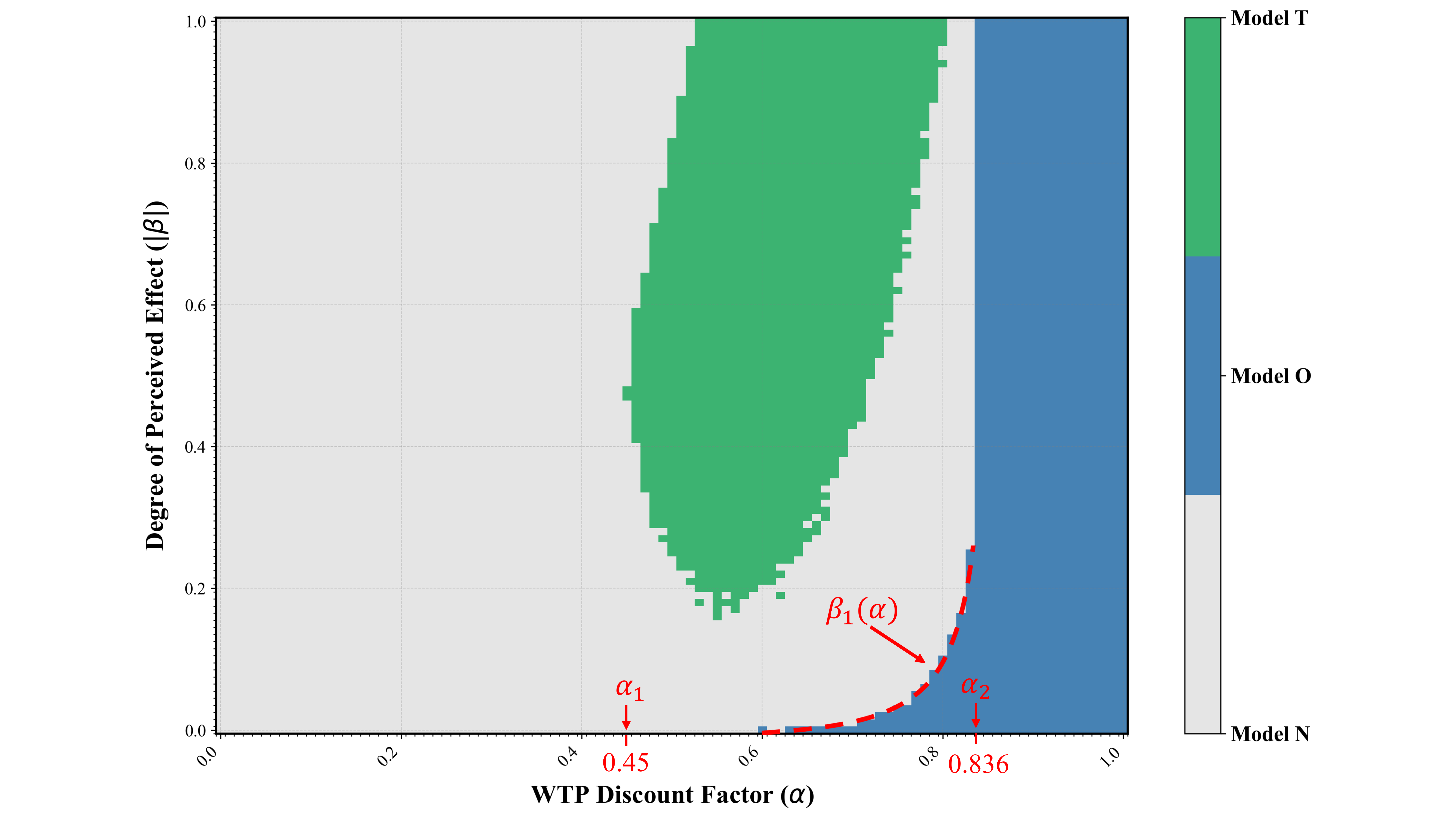}
\caption{Model selection map}
\label{fig:selection map H10k h100}
\end{figure}

Fig.~\ref{fig:selection map H10k h100} illustrates the optimal selection of remanufacturing business models for all consumer perception conditions. It indicates that the consumer's perceived value of the remanufactured product, i.e., $\alpha$, is the main driver for business model selection. The degree of assimilation and contrast effects, i.e., $|\beta|$, becomes decisive only when $\alpha$ is moderate ($\alpha_1\leq \alpha \leq \alpha_2$), especially influencing the selection of TPR-authorized remanufacturing. In this figure, $\beta_1(\alpha)$ represents the approximate boundary between models N and O. It is derived by approximating the expected profit functions for both models by replacing $F(q^*-1)$ with $F(q^*)$, and comparing their resulting approximate profits. See Appendix for derivation details.

% In addition, the observation reveals a distinct non-monotonic transition with respect to $\alpha$ for $\vert \beta \vert >0.2$: as $\alpha$ increases, the selection shifts from Model N to Model T, then reverts back to Model N, and finally transitions to Model O. This pattern is driven by the underlying profit dynamics. Under Model T, the OEM earns from the sales of new products and authorization fees from the TPR. As $\alpha$ increases, remanufactured products attract more customers with low preferences, which initially increases the OEM's revenue through increased authorization fees. However, further increases in $\alpha$ intensify competition between new and remanufactured products, as consumers perceive their values to be increasingly similar. This intensifies price competition, reducing the margin of new products; eventually, this loss in margin exceeds the additional revenue from authorization, causing the OEM’s profit under Model T to decline.
\begin{figure}[ht!]
\centering
\begin{tikzpicture}[scale=0.7,
    box/.style={
        rectangle, 
        draw=black,
        thick,
        minimum width=2cm, 
        minimum height=1cm,
        fill=white
    },
    outcome/.style={
        rectangle, 
        draw=black,
        thick,
        text width=3cm, 
        minimum height=2cm, 
        align=center,
        font=\bfseries
    }
]

% Define colors
\definecolor{n_grey}{HTML}{E5E5E5}
\definecolor{t_green}{HTML}{3CB371}
\definecolor{o_blue}{HTML}{4682B4}

% decision nodes 
\node[box, rounded corners=5pt, fill=white] (alpha) at (0,0) {\textbf{Assess $\alpha$}};
\node[box, rounded corners=5pt, fill=white] (beta) at (6,0) {\textbf{Assess $|\beta|$}};
\node[box, rounded corners=5pt, align=center, fill=white] (contract) at (10,2) 
    {\textbf{\small Contract Design}\\ \small$(h,H)$};

% Outcome nodes
\node[outcome, fill=n_grey] (model_n) at (16,4) {\scriptsize Model N:\\No Remanufacturing};
\node[outcome, fill=t_green] (model_t) at (16,0) {\scriptsize Model T:\\TPR-authorized\\Remanufacturing};
\node[outcome, fill=o_blue] (model_o) at (16,-4) {\scriptsize Model O:\\OEM In-house\\Remanufacturing};

% curved paths
\draw[->, thick] (alpha) to[out=45, in=180] node[above, sloped, font=\footnotesize\bfseries, text width=3cm, align=center] 
    {Low\\ $(0 \leq \alpha \leq \alpha_1)$} (model_n);
\draw[->, thick] (alpha) -- node[ font=\footnotesize\bfseries, text width=3cm, align=center] 
    {Moderate\\$(\alpha_1 \leq \alpha \leq \alpha_2)$} (beta);
\draw[->, thick] (alpha) to[out=-45, in=180] node[below, sloped, font=\footnotesize\bfseries, text width=3cm, align=center] 
    {High\\ $(\alpha_2 \leq \alpha \leq 1)$} (model_o);

\draw[->, thick] (beta) to[out=60, in=180] node[above, sloped, font=\scriptsize\bfseries, yshift = 5pt, xshift = 2pt] 
    {Moderate/high} (contract);
\draw[->, thick] (beta) to[out=-60, in=180] node[above, sloped, font=\scriptsize\bfseries, text width=3cm, align=center] 
    {Low\\ $(0 \leq |\beta| < \beta_1(\alpha))$} (model_o);

% dashed lines
\draw[->, dashed, thick] (contract) to[out=55, in=180] (model_n);
\draw[->, dashed, thick] (contract) to [out=-45, in=180] (model_t);

\node[font=\scriptsize, align=left] at ($(contract.east)+(2.1,0)$) {Depends on contract};

\end{tikzpicture}
\caption{Hierarchical decision roadmap} \label{fig:roadmap}
\end{figure}
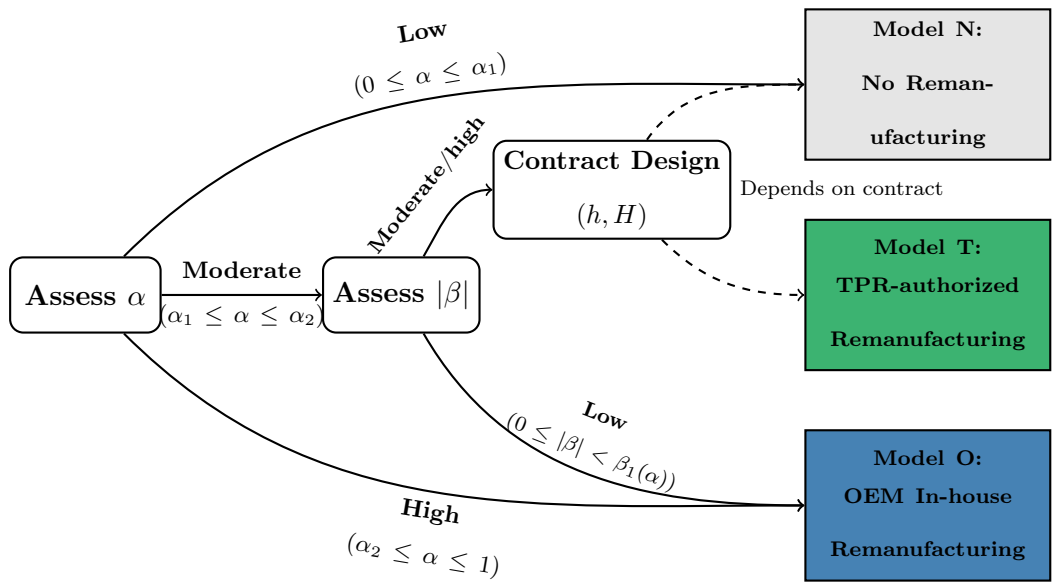

Building on our numerical analysis, we propose a hierarchical decision roadmap in Fig.~\ref{fig:roadmap}, which will guide the OEM through the model selection process based on consumer conditions. This roadmap suggests that the OEM can streamline their market research by first assessing consumer attitudes towards the remanufactured product ($\alpha$) in their target market. If the attitude is highly negative or positive, the optimal business decision is clear, that is, no remanufacturing or OEM in-house remanufacturing, respectively. For moderate attitude, the OEM should then assess the magnitude of assimilation and contrast effects ($\vert \beta \vert$) to determine the most profitable model. 

To assess the practical relevance of our findings, we review empirical studies, indicating that typical values of $\alpha$ range from $0.4$ to $0.9$ across different product categories \citep{guide2010potential,subramanian2012key,abbey2017role}. However, empirical research on assimilation and contrast effects ($\beta$) remains limited. The only study on customer goods, such as MP3 players and printers, reports a value for $\beta$ of around 8\%, while corresponding values for other product categories have not been documented yet. Based on these findings, we identify a realistic parameter zone of $[0.4,0.9] \times [0, 0.3]$. As illustrated in Fig.~\ref{fig:selection map H10k h100}, the practical range spans all three remanufacturing business models, thus validating our contribution as a guide for OEMs.

\subsection{Market dynamics}
\begin{figure}[ht!]
\centering
\begin{subfigure}[b]{.3\linewidth}
\includegraphics[width=\linewidth]{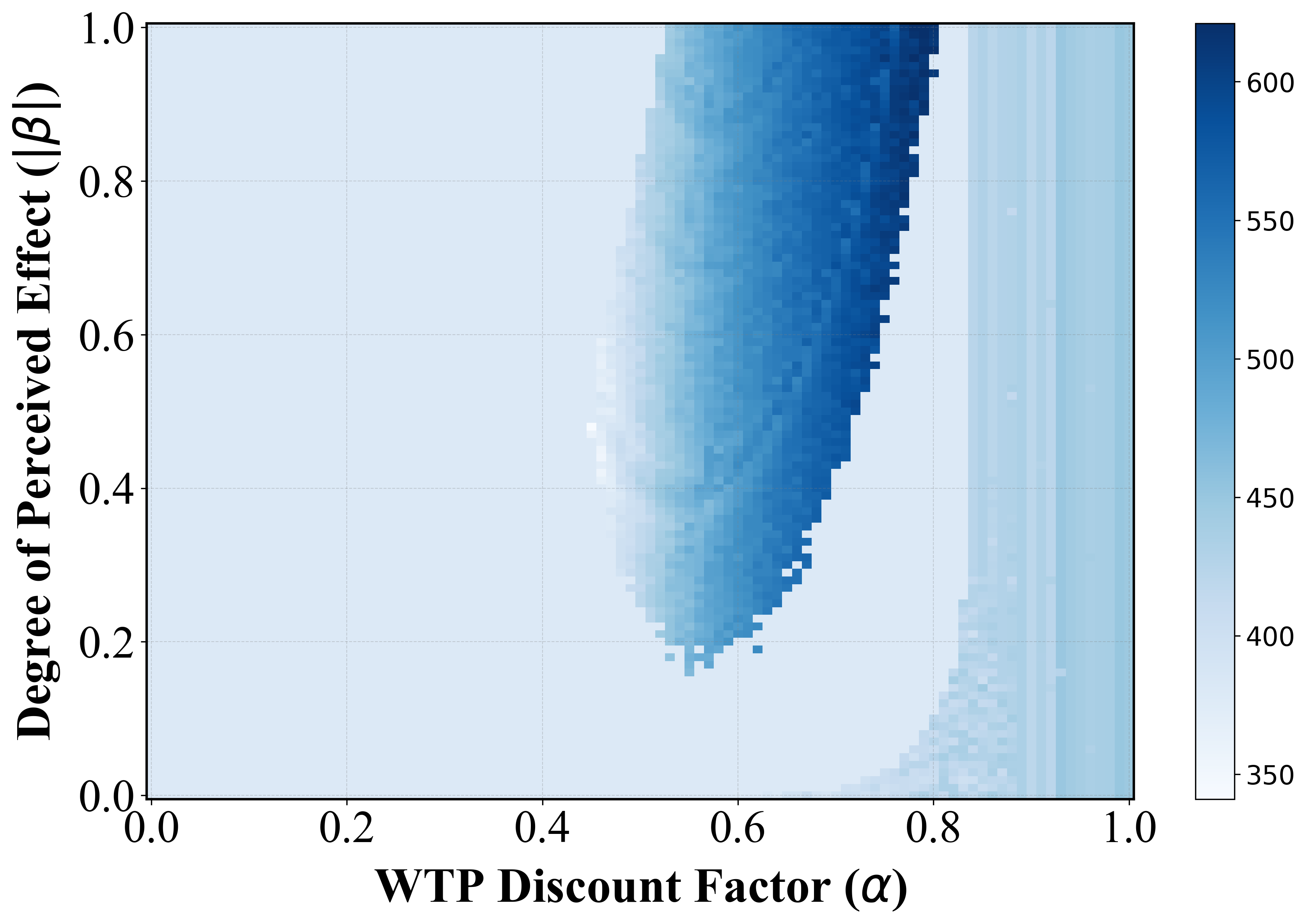}
\caption{Total product quantities}\label{fig:total q}
\end{subfigure}%
\hfill
\begin{subfigure}[b]{.3\linewidth}
\includegraphics[width=\linewidth]{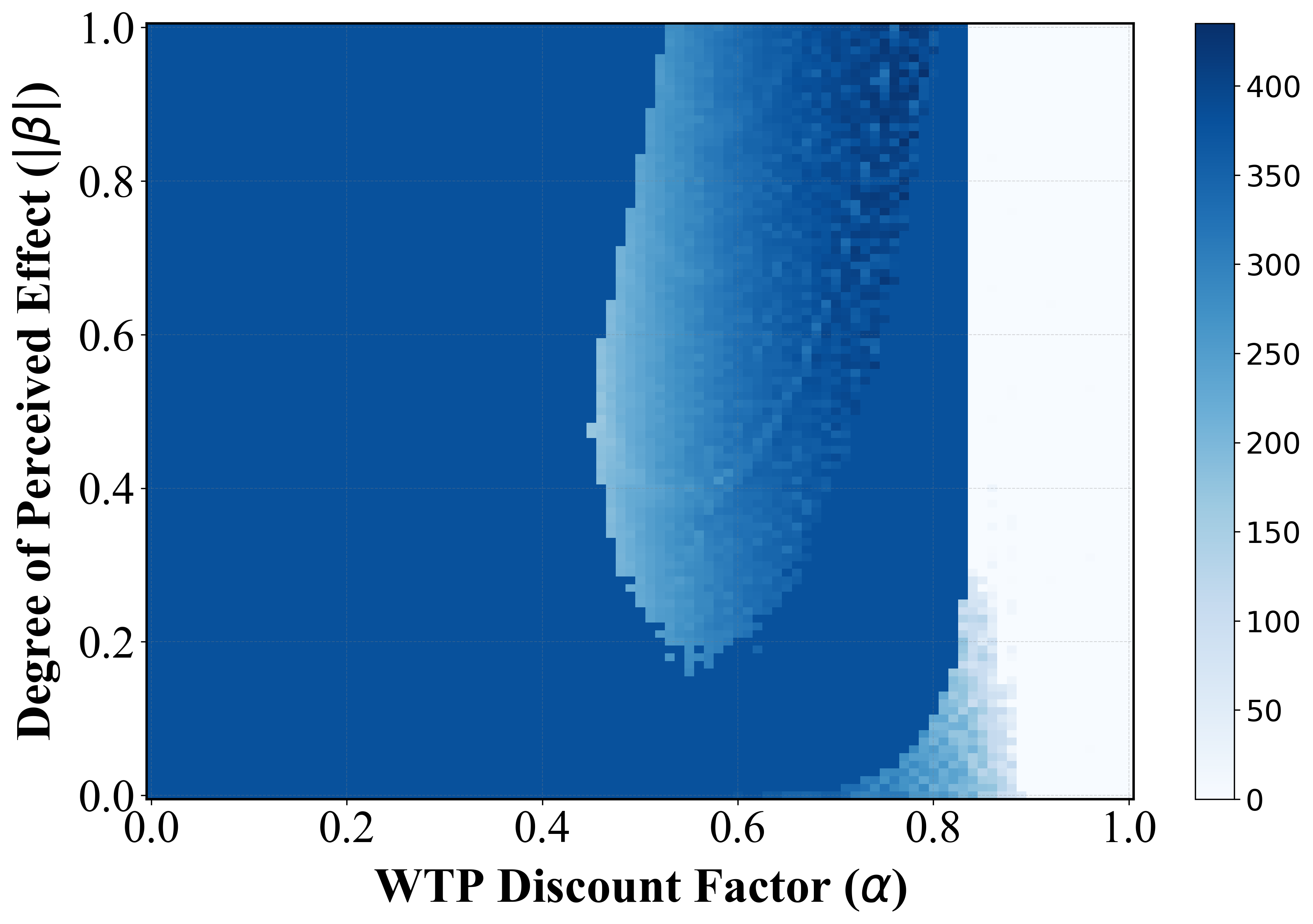}
\caption{Quantity of new products}\label{fig:qn}
\end{subfigure}
\hfill
\begin{subfigure}[b]{.3\linewidth}
\includegraphics[width=\linewidth]{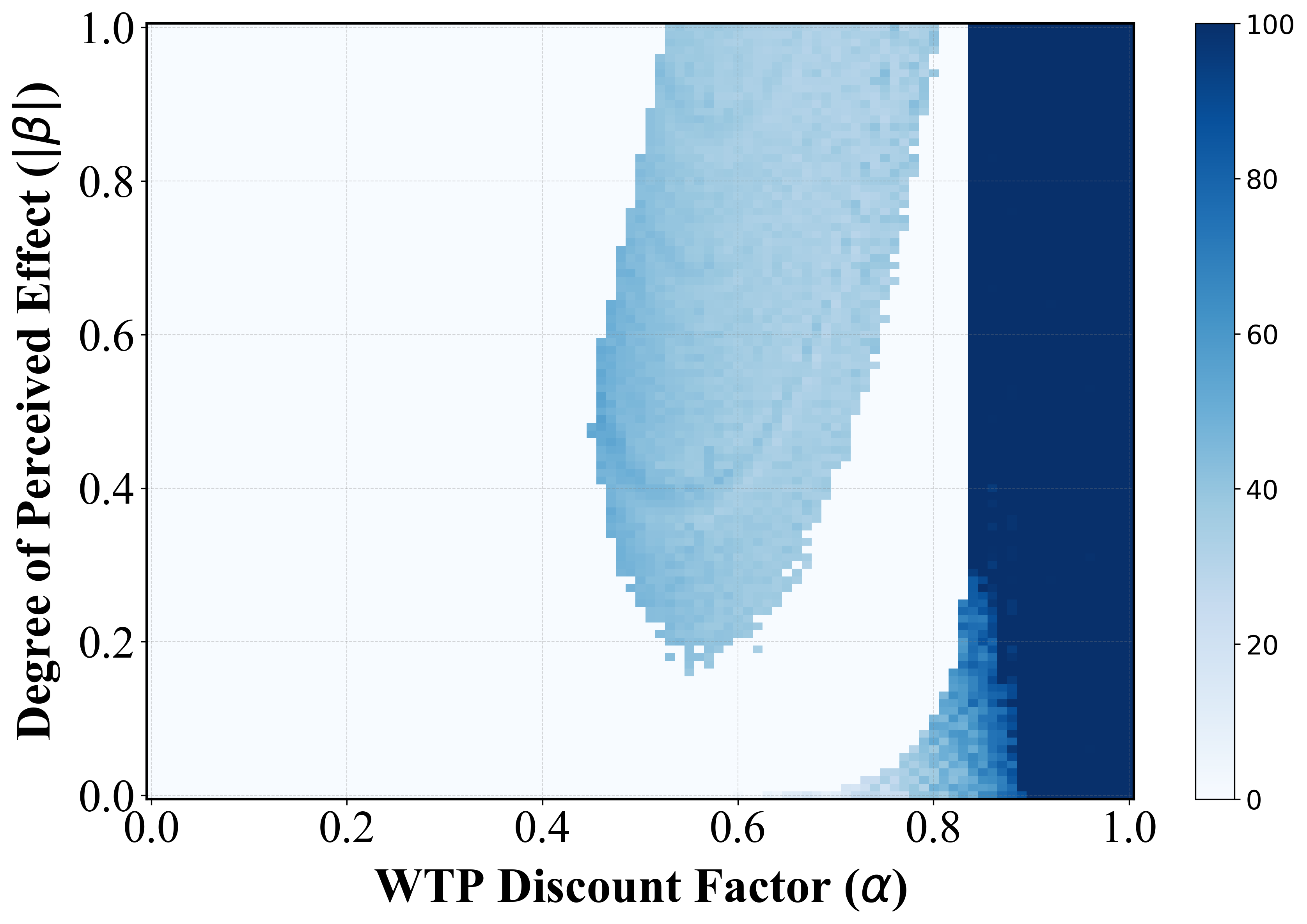}
\caption{Proportion of reman. products}\label{fig:qr percentage}
\end{subfigure}

\caption{Market dynamic under optimal remanufacturing model}\label{fig:market dynamic}
\end{figure}

In this section, we analyze how the selection of remanufacturing business models influences subsequent market outcomes. To capture the market dynamics, we conduct numerical studies to evaluate three key metrics: (1) the total quantities of new and remanufactured products ($q_n+q_r$), (2) the quantity of new products ($q_n$), and (3) the proportion of remanufactured products relative to the total market ($\frac{q_r}{q_n+q_r}\times 100\%$), see Fig.~\ref{fig:market dynamic}. 

Fig.~\ref{fig:total q} shows that the optimal model typically offers benefits in market expansion. Due to the differentiation in value between new and remanufactured products, more customers with lower preferences participate in the market. In particular, the total quantities can increase to 63.9\% when the
TPR-authorized remanufacturing model is optimal.
Fig.~\ref{fig:qn} reveals that, contrary to the common concern of managers about cannibalization of new products \citep{guide2010potential, atasu2010so}, remanufacturing can increase sales of new products by up to 14.8\%, when TPR-authorized remanufacturing is optimal. The reason is that the presence of TPR-remanufactured products can significantly enhance the perceived value of new products and attract more customers when the contrast effect is strong.
Fig.~\ref{fig:qr percentage} highlights a notable market outcome. When customers highly value remanufactured products (very high $\alpha$), the market can become dominated by OEM-remanufactured products. 
This market scenario aligns with findings of \citet{atasu2008remanufacturing}, in which environmentally conscious consumers buy exclusively remanufactured products. Moreover, this outcome challenges the common assumption that new and remanufactured products must necessarily coexist in the market \citep{atasu2008remanufacturing,xiong2013don,zou2016third}.

\subsection{Environmental impact}
We examine how OEM's optimal remanufacturing model influences environmental outcomes. 
Based on the Life Cycle Assessment (LCA) of \citet{orsdemir2014competitive,zou2016third} and \citet{jin2023right}, we evaluate the environmental impacts in three phases: (re)production, consumer consumption and disposal. Each phase has its own unit environmental impact, denoted as $e_p$ (or $e_r$) for production or remanufacturing, $e_c$ for consumption, and $e_d$ for disposal. Restricting $e_r<e_p$ reflects the environmental benefit of remanufacturing process, as it incorporates reused materials compared to producing a new item.
The overall environmental impact depends on the quantities of new and remanufactured products, i.e., $q_n$ and $q_r$, and their respective sales, i.e., $S_n,S_r$, which are related to the consumer consumption phase. Consequently, the total environmental impact (EI) is defined as follows:
\begin{align*}
    EI & := e_p  q_n + e_r  q_r + e_c  \mathbb{E}(S_n + S_r) + e_d (q_n+q_r)\\
         & = \underbrace{(e_p+e_d)}_{\gamma_n} q_n + e_c \mathbb{E}(S_n) + \underbrace{(e_r+e_d)}_{\gamma_r} q_r + e_c  \mathbb{E}(S_r).
\end{align*}
For convenience, we denote $\gamma_n:= e_p+e_d$ as the combined impact of production and disposal of a new product, and $\gamma_r:= e_r+e_d$ as the corresponding value of a remanufactured product.

\begin{figure}[ht!]
\begin{subfigure}{0.45\textwidth}
\centering
\includegraphics[width=0.9\textwidth]{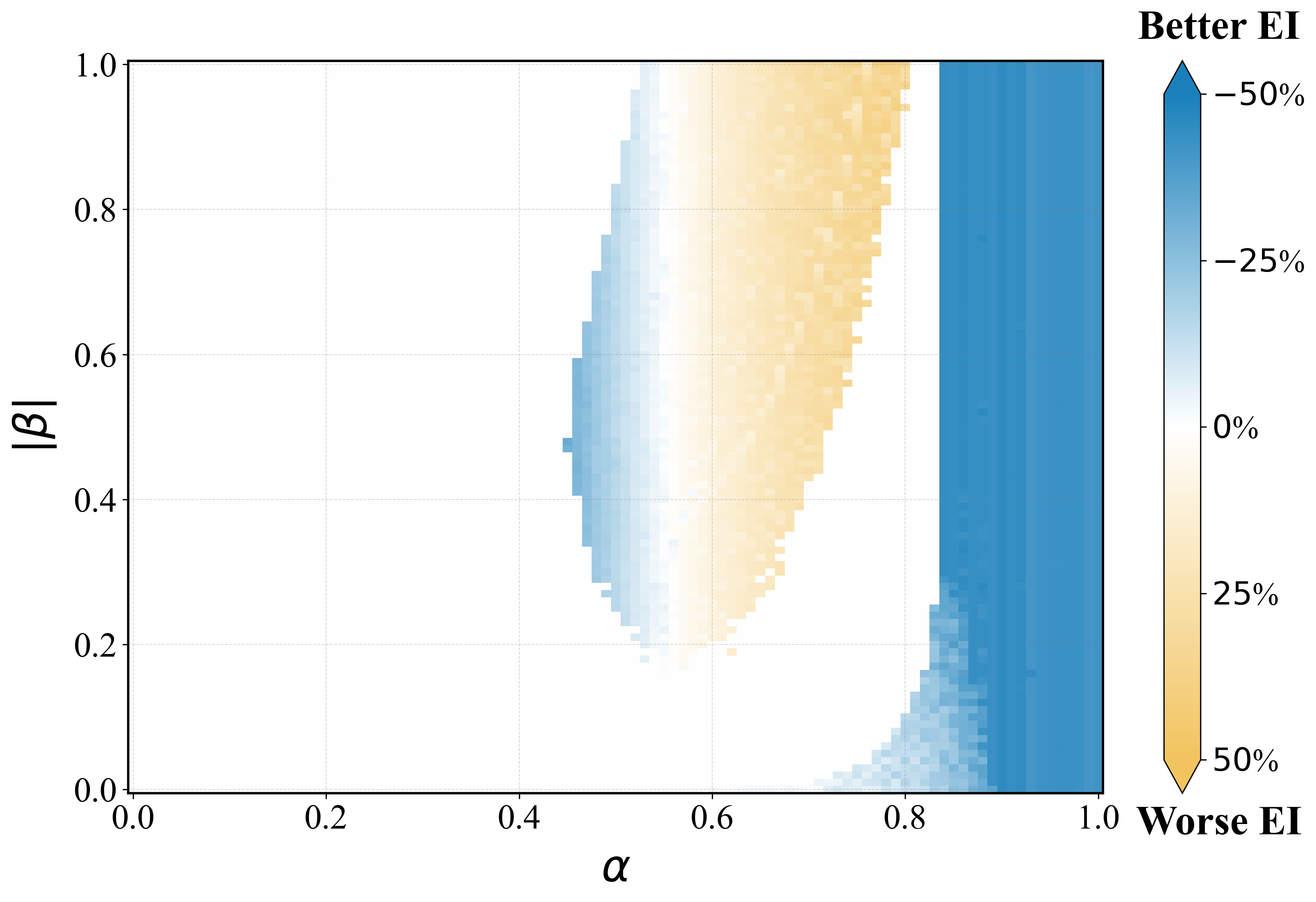}
\caption{Production and disposal phases dominance} 
\label{fig:EI p+d>u}
\end{subfigure}
\hfill
\begin{subfigure}{0.45\textwidth}
\centering
\includegraphics[width=0.9\textwidth]{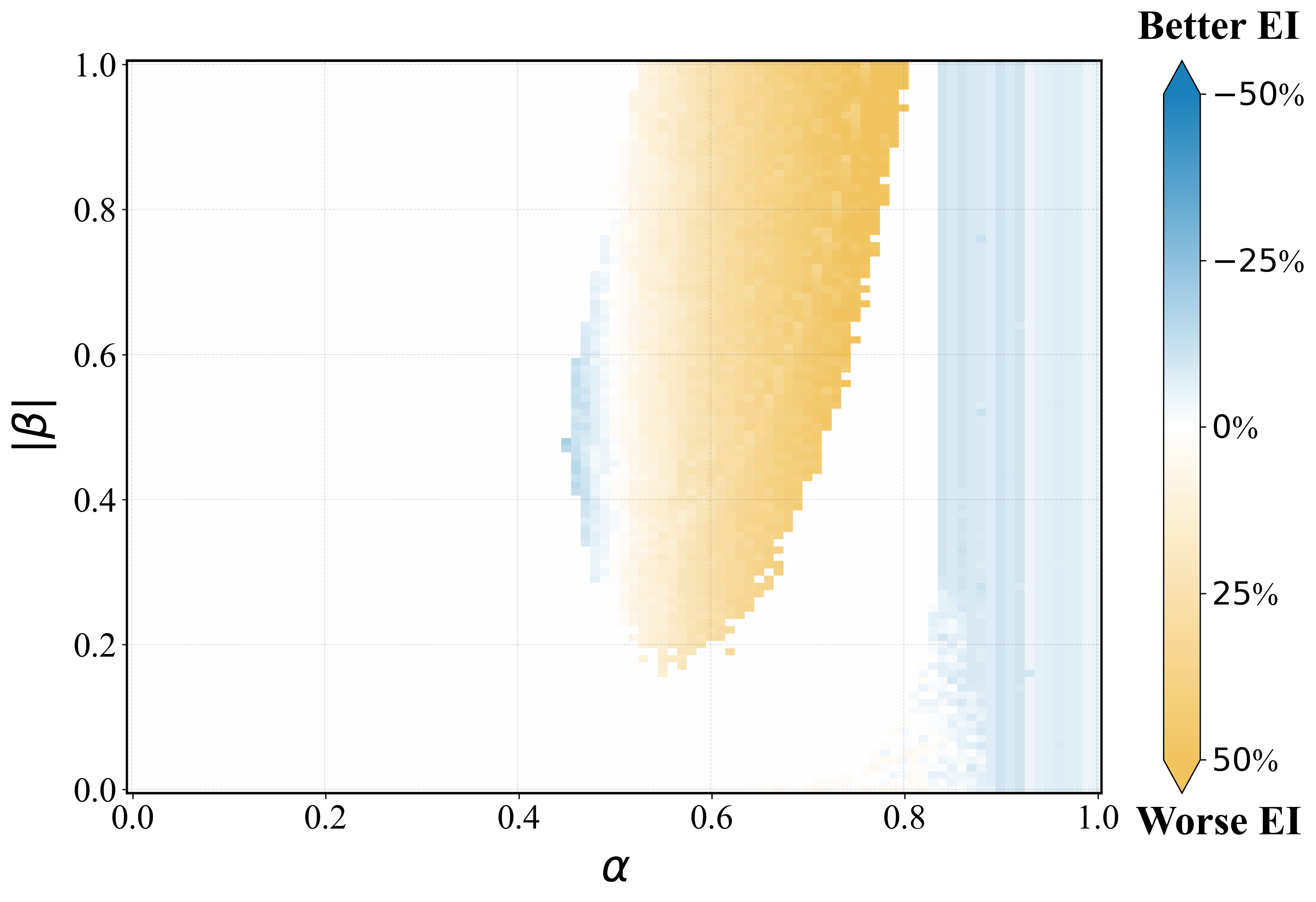}
\caption{Consumption phase dominance} \label{fig:EI p+d<u}
\end{subfigure}
\caption{Environmental impact under optimal remanufacturing model} \label{fig:EI}
\end{figure}

To examine the environmental outcomes across product categories, we evaluate two scenarios based on \citet{zou2016third} and \citet{jin2023right}: (1) production and disposal dominance ($\gamma_n = 7, \gamma_r = 3, e_c=1$), representing, for example, high-end computers; (2) consumption phase dominance ($\gamma_n = 4, \gamma_r = 2, e_c=7$), representing, for example, vehicles and refrigerators. Fig.~\ref{fig:EI} illustrates the environmental outcomes in the two scenarios evaluated. 

By comparing the results of two scenarios, we find that a higher unit impact in the consumption phase leads to a worse overall environmental outcome than higher unit impacts in the production and disposal phases. This highlights that products with consumption-dominant impacts, such as vehicles and refrigerators, face greater challenges in achieving environmental benefits through remanufacturing. Importantly, the findings demonstrate that remanufacturing does not always guarantee environmental advantages. Model O consistently reduces environmental impacts, mainly by substantially decreasing the output of new products and the total quantities of the market. However, Model T tends to increase total environmental impact due to overall market expansion. This effect is most pronounced for products with high impacts in the consumption phase, where market growth can offset resource savings from remanufacturing. Overall, these results highlight a paradox: while remanufacturing is promoted for its environmental benefits, it can, under certain market conditions and especially for usage-intensive products, increase total environmental harm due to induced market expansion and higher aggregate consumption.

\section{Sensitivity analysis}
\label{sec:contract and stochasticity}
In this section, we investigate how the design of the authorization contract and the stochasticity of the market size affect the system profitability and environmental impact.

\subsection{Effect of the authorization contract}
Throughout this study, we assume that the authorization contract adopts a two-part tariff structure, characterized by exogenous parameters $(H,h)$. In this section, we aim to explore the impact of different authorization contracts on the system outcome, rather than to optimize the contract parameters. Specifically, we examine the impact of varying the one-time fee ($H$) and unit fee ($h$) on the system profitability and the environmental impacts. In particular, setting $H=0$ produces a one-part contract, which allows us to examine the impact of different contracts, i.e., one-part and two-part tariff.

\begin{figure}[ht!]
\begin{subfigure}{0.45\textwidth}
\centering
\includegraphics[width=0.9\textwidth]{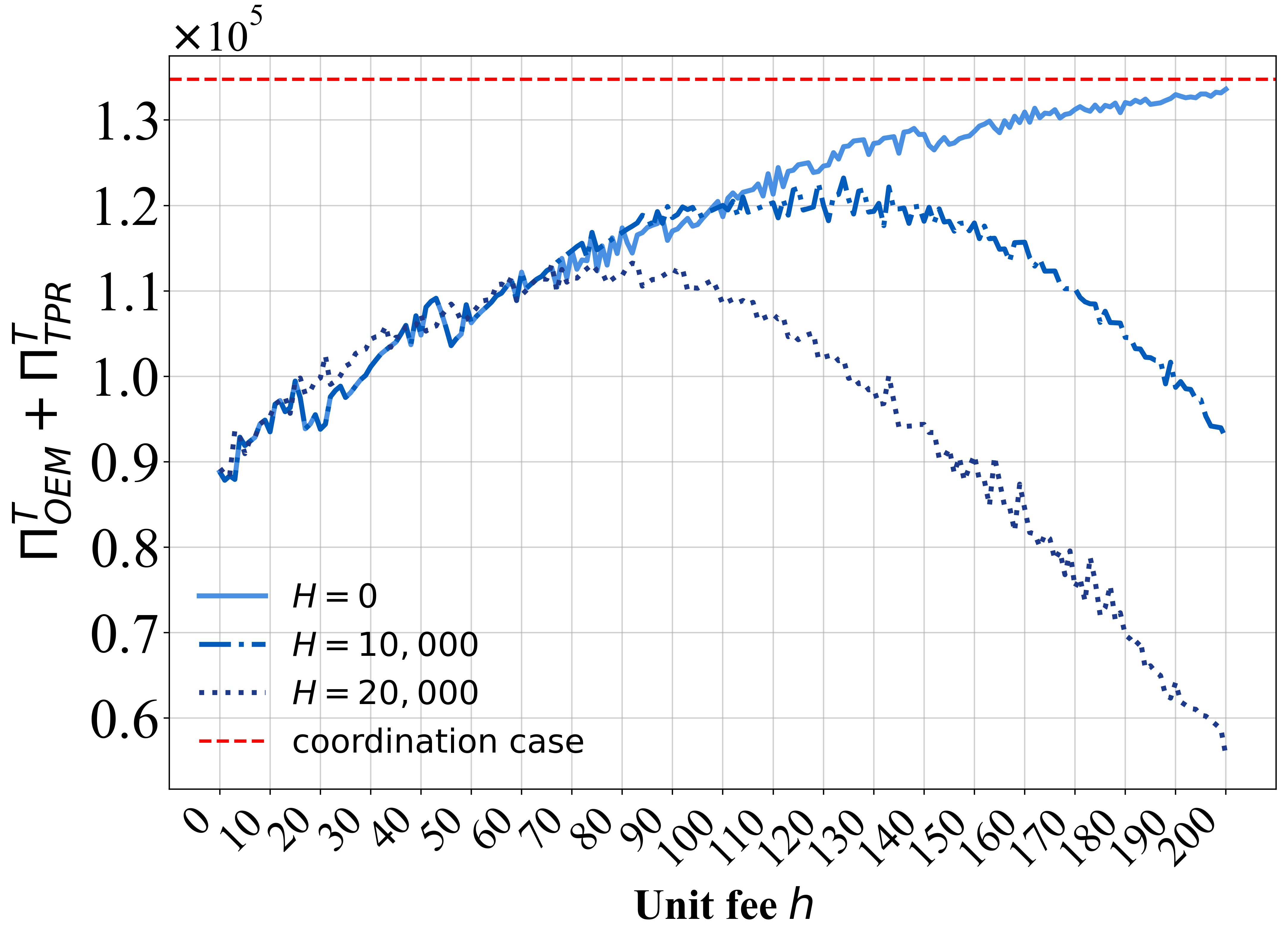}
    \caption{System profitability}
    \label{fig:Hh on profit}
\end{subfigure}
\hfill
\begin{subfigure}{0.45\textwidth}
    \centering
\includegraphics[width=0.9\textwidth]{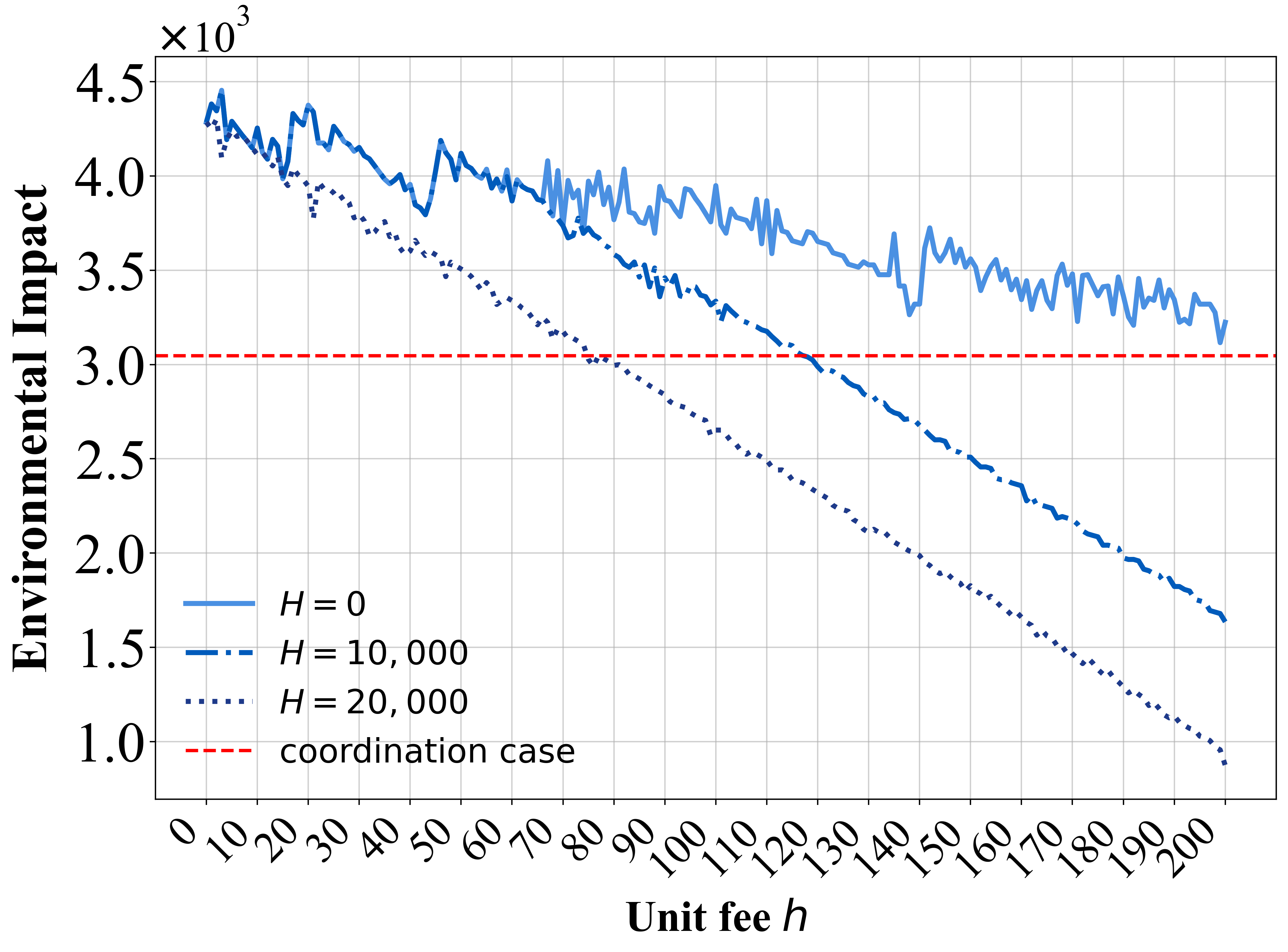}
    \caption{Environmental impact}
    \label{fig:Hh on EI}
\end{subfigure}
\caption{Effect of authorization contract ($\alpha=0.6,\beta^+ =0.3$)}
\label{fig:Hh analysis}
\end{figure}

Fig.~\ref{fig:Hh analysis} shows two metrics: (1) the total profit of OEM and TPR, and (2) the environmental impact, across varying one-time fees and unit fees, i.e., $H=0, 10000, 20000$ and $h\in[0,c]$. We use environmental impact parameters $\gamma_n=7,\gamma_r=3,e_c=1$. Other parameter settings are the same as before. The ``coordination case'' reflects the outcome of optimal decisions in centralized system with contrast effect. For the other setting of environmental impact parameters, i.e., $\gamma_n=4,\gamma_r=2,e_c=7$, the results are structurally the same. 

By comparing outcomes of one-part contracts and two-part contracts, we suggest that, while two-part tariff contracts may not always maximize total system profit, they are feasible to meet more strict environmental requirements compared to one-part contracts. For system profitability, when the unit fee $h$ is low, the results are identical regardless of the value of the one-time fee $H$. As the unit fee increases, profits under two-part contracts decline, while profits under one-part continue to rise and eventually approach those of the centralized system. This contrasts with the conventional view that a two-part tariff contract typically improves coordinate the supply chain due to the additional parameter \citep{corbett2004designing}. In our context, under the two-part contracts, the TPR must generate sufficient remanufactured products to cover the one-time fee, which limits its flexibility in decisions and diminishes overall profit. 
However, this constraint also shifts market share away from new products towards remanufactured products, leading to lower per-unit environmental impacts and decreased total output. As a result, two-part contracts produce lower environmental impacts than one-part contracts, making them particularly suitable when strict environmental restrictions must be met. 

\subsection{Effect of stochastic market size}
We model the stochasticity of the market size and assume a Poisson distribution, i.e., $N\sim \text{Poisson}(\lambda)$. To evaluate the effects of stochasticity on system outcomes, we first analyze the case of a constant market size, i.e., $N=\lambda$, and derive its optimal remanufacturing decisions, denoted by $d^{cons}$. Next, we apply these deterministic decisions to the stochastic scenario and compute the resulting outcomes, denoted by $\Pi(d^{cons})$ and $EI(d^{cons})$. We also derive the optimal decisions under the stochastic case and compute the resulting outcomes, denoted by $\Pi(d^{stoc})$ and $EI(d^{stoc})$. 
\begin{figure}[ht!]
\begin{subfigure}{0.45\textwidth}
\centering
\includegraphics[width=0.9\textwidth]{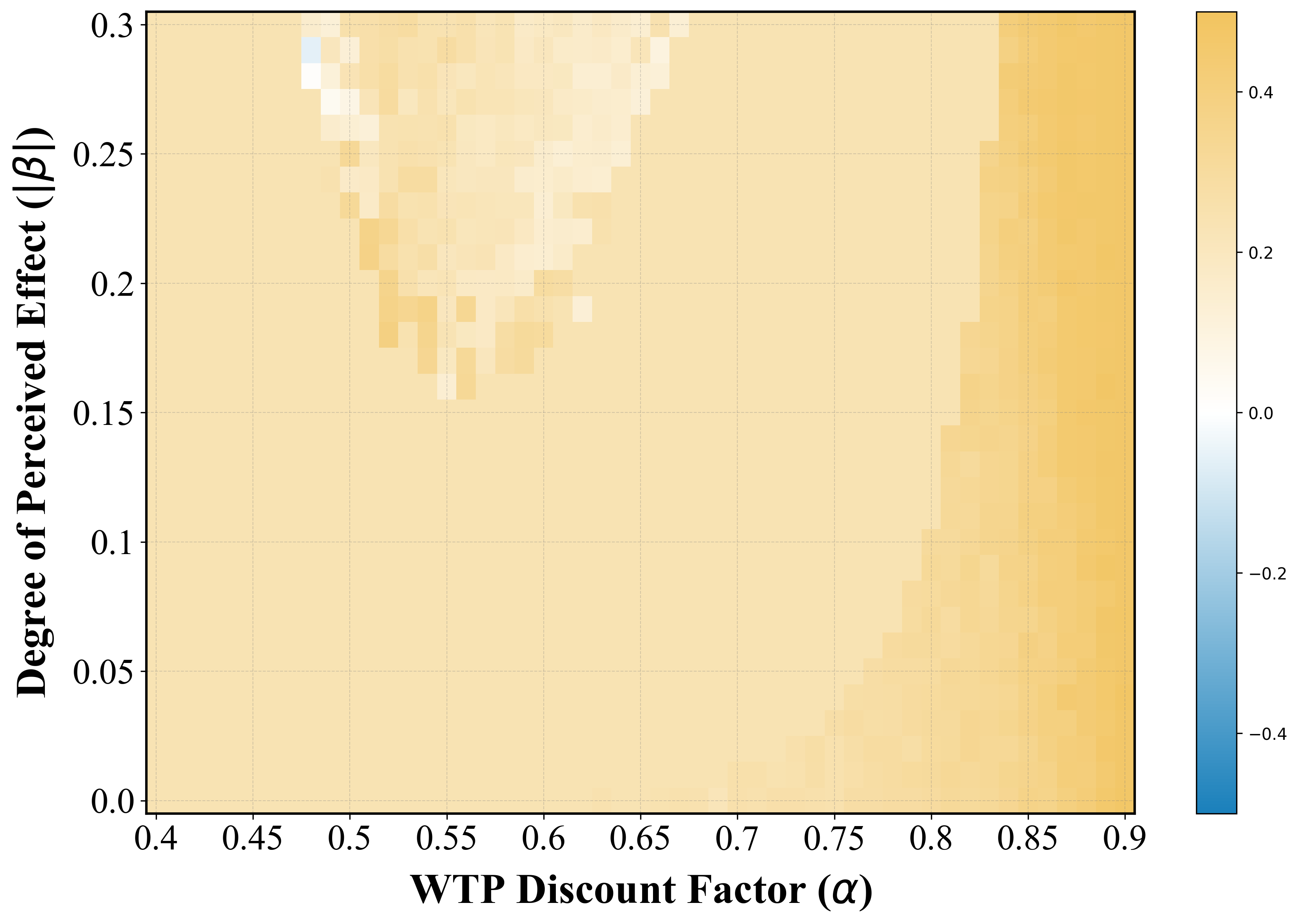}
\caption{Relative differences in profitability}
\label{fig:stochastic market size on profit}
\end{subfigure}
\hfill
\begin{subfigure}{0.45\textwidth}
    \centering
\includegraphics[width=0.9\textwidth]{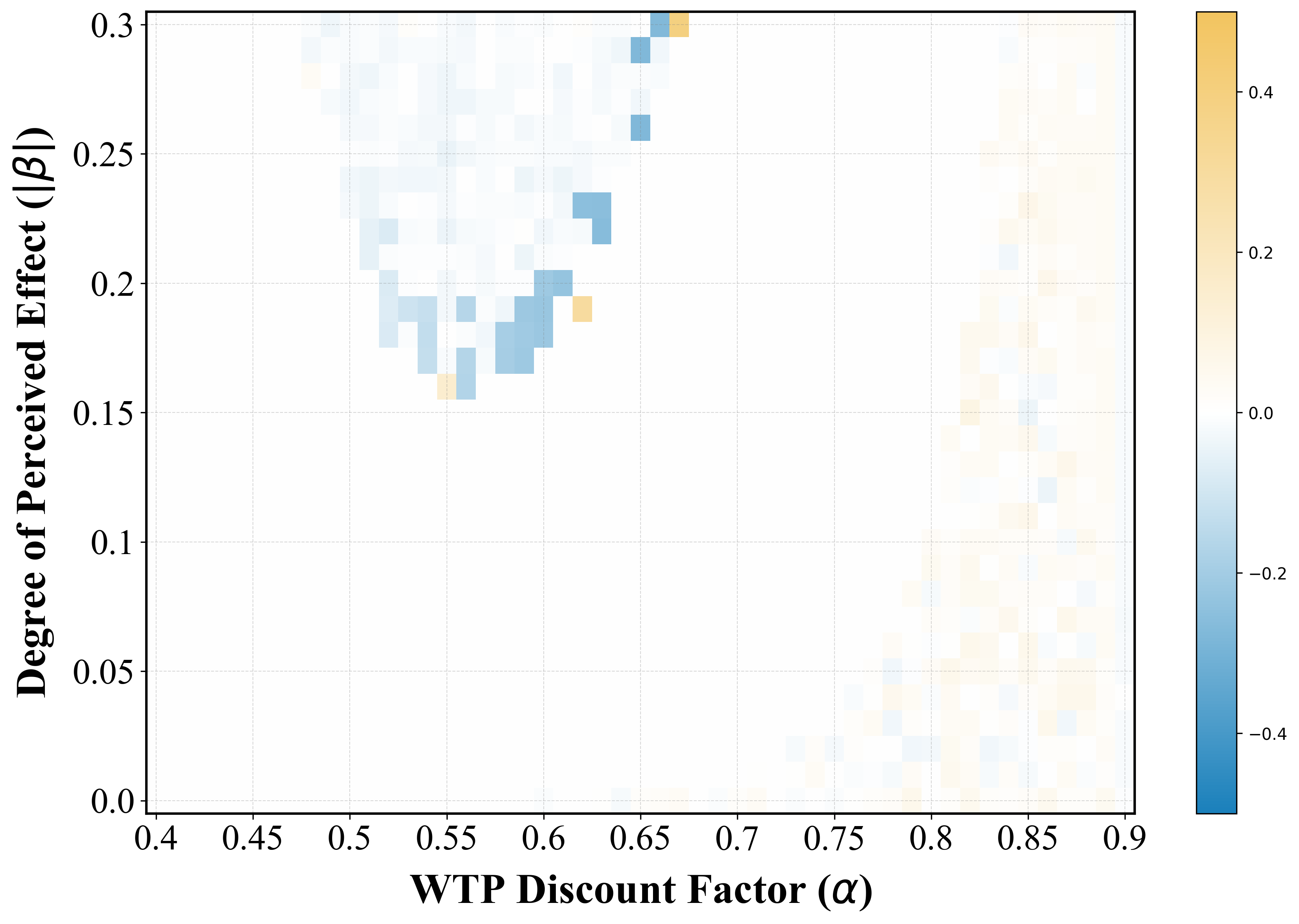}
\caption{Relative differences in environmental impact}
\label{fig:stochastic market size on EI}
\end{subfigure}
\caption{Effect of stochastic market size}
\label{fig:stochastic market size analysis}
\end{figure}

Fig.~\ref{fig:stochastic market size analysis} shows the relative change of the outcomes from constant to stochastic decisions, i.e., $\frac{\Pi(d^{stoc})-\Pi(d^{cons})}{\Pi(d^{cons})}$ for profitability and $\frac{EI(d^{stoc})-EI(d^{cons})}{EI(d^{cons})}$ for environmental impact, in the realistic parameter zone ($[0.4,0.9]\times[0,0.3]$). Here, the environmental impact parameters are $\gamma_n=4,\gamma_r=2,e_c=7$. 
Fig.~\ref{fig:stochastic market size on profit} indicates that accounting for the stochasticity of market size consistently improves system profitability, regardless of consumer perception conditions. The overall system profitability may increase by up to 44.8\%. Interestingly, the impact on the environmental impact varies. In most cases, the TPR-authorized model (top-left colored area in Fig.~\ref{fig:stochastic market size on EI}) reduces environmental impact, while the OEM in-house model (right colored area in Fig.~\ref{fig:stochastic market size on EI}) can result in a slight increase, typically less than 5\%. This difference stems from overproduction in the TPR-authorized model and underproduction in the OEM in-house model under the constant case. Note that minor numerical deviations near selection boundaries are attributed to the discrete nature of decision variables and are therefore ignored when analyzing the overall pattern. In summary, these findings suggest that incorporating the stochasticity of market size into remanufacturing optimization can deliver economic improvement without substantially compromising environmental impacts.

\section{Conclusions} \label{sec:conclusion}

This study investigates how consumer perceptions determine the OEM's selection of the most profitable remanufacturing business models. Consumer perceptions consist of the WTP discount factor for remanufactured products, i.e., $\alpha$, and the degree of assimilation and contrast effects, i.e., $\vert \beta \vert$. We analyze three alternative models: no remanufacturing (Model N), OEM in-house remanufacturing (Model O) and TPR-authorized remanufacturing (Model T) with a stochastic market size and a two-part tariff authorization contract. The numerical results provide managerial insights for model selection. 
The OEM should first evaluate consumer attitudes towards remanufactured products in their target market. If attitudes are unfavorable, forgoing remanufacturing is optimal for the OEM; if attitudes are highly favorable, OEM in-house remanufacturing is most profitable. For intermediate attitudes, the OEM should assess the magnitude of assimilation and contrast effects. If the magnitude is not very low, TPR-authorized remanufacturing generally emerges as the optimal model. 

Beyond OEM's profitability, we examine market and environmental implications under the optimal remanufacturing strategy. From a market perspective, optimal remanufacturing models generally expand the total market sales by attracting new customer segments. Contrary to traditional concerns of cannibalization, TPR-authorized remanufacturing can even increase sales of new products through the contrast effect, representing a ``win-win-lose'' outcome: higher profitability, substantial market expansion, but increased environmental burden. 
In contrast, when the OEM in-house remanufacturing model is optimal, a ``win-win-win'' outcome arises, characterized by higher profits, moderate market growth, and lower environmental impact. In this case, remanufactured products may dominate the market. This finding is observed in industry practice. For example, Apple’s exclusive remanufactured sales of certain products, such as iPhone 15 and MacBook M3, illustrate that remanufacturing dominance already occurs in practice. These results highlight that although remanufacturing aims to promote sustainability, it can eventually increase the total environmental burden due to the joint decisions on price and quantity.

Our study further examines how the design of the authorization contract and the stochasticity of market size affect the system's profitability and environmental impact. The numerical results show that, while the one-part structure of authorization contract generates higher profits in most cases, it produces a higher environmental impact compared to the two-part structure. In contrast, two-part tariff contracts are more suitable for balancing profitability and sustainability and are more feasible to satisfy strict environmental impact caps. In addition, although demand uncertainty has been studied by heterogeneous consumer preferences in prior research, our findings highlight that accounting for the stochasticity of market size, an often overlooked factor, can further improve system profitability without significantly compromising sustainability.

We conclude by highlighting potential directions for future research. This study employs a single-period model that examines competition between new and remanufactured products. As a result, it does not account for the dynamic nature of remanufacturing decisions over time, which may be particularly relevant in industries characterized by long product lifecycles or active secondary markets. Additionally, for the sake of analytical tractability, we assume linear collection costs; however, real-world collection processes may involve economies of scale that warrant further investigation.

\section*{Acknowledgements}
This research was conducted as part of project \textit{Transitioning to a Circular Business Ecosystem} (LINCIT)\, supported by the TKI Dinalog. 
% \textcolor{red}{This activity is co-funded by the Temporary research subsidy scheme Top Sector Logistics 2022-2026 of the Ministry of Infrastructure and Water Management.} 
This work used the Dutch national e-infrastructure with the support of the SURF Cooperative using grant no. EINF-15313.

\bibliographystyle{apalike-ejor}
\linespread{1.5}
\bibliography{ref.bib}

@article{agrawal2015remanufacturing,
  title={Remanufacturing, third-party competition, and consumers' perceived value of new products},
  author={Agrawal, Vishal and Atasu, Atalay and Van Ittersum, Koert},
  journal={Management Science},
  volume={61},
  number={1},
  pages={60--72},
  year={2015},
  publisher={INFORMS}
}

@article{zou2016third,
  title={Third-party remanufacturing mode selection: Outsourcing or authorization?},
  author={Zou, Zongbao and Wang, Jianjun and Deng, Guishi and Chen, Haozhe},
  journal={Transportation Research Part E: Logistics and Transportation Review},
  volume={87},
  pages={1--19},
  year={2016},
  publisher={Elsevier}
}

@article{guide2010potential,
  title={The potential for cannibalization of new products sales by remanufactured products},
  author={Guide Jr, V Daniel R and Li, Jiayi},
  journal={Decision Sciences},
  volume={41},
  number={3},
  pages={547--572},
  year={2010},
  publisher={Wiley Online Library}
}

@article{wu2020competitive,
  title={Competitive remanufacturing and pricing strategy with contrast effect and assimilation effect},
  author={Wu, Lingli and Liu, Liwen and Wang, Zongjun},
  journal={Journal of Cleaner Production},
  volume={257},
  pages={120333},
  year={2020},
  publisher={Elsevier}
}

@article{donohue2020behavioral,
  title={Behavioral operations: Past, present, and future},
  author={Donohue, Karen and {\"O}zer, {\"O}zalp and Zheng, Yanchong},
  journal={Manufacturing \& Service Operations Management},
  volume={22},
  number={1},
  pages={191--202},
  year={2020},
  publisher={INFORMS}
}

@article{subramanian2012key,
  title={Key factors in the market for remanufactured products},
  author={Subramanian, Ravi and Subramanyam, Ramanath},
  journal={Manufacturing \& Service Operations Management},
  volume={14},
  number={2},
  pages={315--326},
  year={2012},
  publisher={INFORMS}
}

@article{huang2024pride,
  title={Pride or guilt? Impacts of consumers’ socially influenced recycling behaviors on closed-loop supply chains},
  author={Huang, Wenjie and Nguyen, Jason and Tseng, Chung-Li and Chen, Wenlin and Kirshner, Samuel N},
  journal={Manufacturing \& Service Operations Management},
  volume={26},
  number={6},
  pages={2067--2084},
  year={2024},
  publisher={INFORMS}
}

@article{abbey2017role,
  title={The role of perceived quality risk in pricing remanufactured products},
  author={Abbey, James D and Kleber, Rainer and Souza, Gilvan C and Voigt, Guido},
  journal={Production and Operations Management},
  volume={26},
  number={1},
  pages={100--115},
  year={2017},
  publisher={SAGE Publications Sage CA: Los Angeles, CA}
}

@article{abbey2015remanufactured,
  title={Remanufactured products in closed-loop supply chains for consumer goods},
  author={Abbey, James D and Meloy, Margaret G and Guide Jr, V Daniel R and Atalay, Selin},
  journal={Production and Operations Management},
  volume={24},
  number={3},
  pages={488--503},
  year={2015},
  publisher={SAGE Publications Sage CA: Los Angeles, CA}
}

@article{ovchinnikov2011revenue,
  title={Revenue and cost management for remanufactured products},
  author={Ovchinnikov, Anton},
  journal={Production and Operations Management},
  volume={20},
  number={6},
  pages={824--840},
  year={2011},
  publisher={SAGE Publications Sage CA: Los Angeles, CA}
}

@article{corbett2004designing,
  title={Designing supply contracts: Contract type and information asymmetry},
  author={Corbett, Charles J and Zhou, Deming and Tang, Christopher S},
  journal={Management Science},
  volume={50},
  number={4},
  pages={550--559},
  year={2004},
  publisher={INFORMS}
}

@article{shi2020distribution,
  title={Distribution channel choice and divisional conflict in remanufacturing operations},
  author={Shi, Tianqin and Chhajed, Dilip and Wan, Zhixi and Liu, Yunchuan},
  journal={Production and Operations Management},
  volume={29},
  number={7},
  pages={1702--1719},
  year={2020},
  publisher={SAGE Publications Sage CA: Los Angeles, CA}
}

@article{liu2022rent,
  title={Rent, sell, and remanufacture: The manufacturer's choice when remanufacturing can be outsourced},
  author={Liu, Jian and Mantin, Benny and Song, Xuefeng},
  journal={European Journal of Operational Research},
  volume={303},
  number={1},
  pages={184--200},
  year={2022},
  publisher={Elsevier}
}

@article{abbey2024closed,
  title={Closed-loop supply chains: a strategic overview},
  author={Abbey, James D and Guide Jr, V Daniel R and Sun, Xichen},
  journal={Sustainable Supply Chains: A Research-based Textbook on Operations and Strategy},
  pages={355--378},
  year={2024},
  publisher={Springer}
}

@article{fang2020third,
  title={Is third-party remanufacturing necessarily harmful to the original equipment manufacturer?},
  author={Fang, Chang and You, Zhuangzhuang and Yang, Yudou and Chen, Duomei and Mukhopadhyay, Samar},
  journal={Annals of Operations Research},
  volume={291},
  pages={317--338},
  year={2020},
  publisher={Springer}
}

@phdthesis{agrawal2010essays,
  title={Essays on sustainable operations},
  author={Agrawal, Vishal},
  year={2010},
  school={Georgia Institute of Technology}
}

@article{lin2024house,
  title={In-house or outsourcing? The impact of remanufacturing strategies on the dynamics of component remanufacturing systems under lifecycle demand and returns},
  author={Lin, Junyi and Naim, Mohamed M and Tang, Ou},
  journal={European Journal of Operational Research},
  volume={315},
  number={3},
  pages={965--979},
  year={2024},
  publisher={Elsevier}
}

@article{feng2021environmentally,
  title={Environmentally responsible closed-loop supply chain models with outsourcing and authorization options},
  author={Feng, Zhangwei and Xiao, Tiaojun and Robb, David J},
  journal={Journal of Cleaner Production},
  volume={278},
  pages={123791},
  year={2021},
  publisher={Elsevier}
}

@article{orsdemir2014competitive,
  title={Competitive quality choice and remanufacturing},
  author={{\"O}rsdemir, Adem and Kemahl{\i}o{\u{g}}lu-Ziya, Eda and Parlakt{\"u}rk, Ali K},
  journal={Production and Operations Management},
  volume={23},
  number={1},
  pages={48--64},
  year={2014},
  publisher={SAGE Publications Sage CA: Los Angeles, CA}
}

@article{wang2017design,
  title={Design of the reverse channel for remanufacturing: must profit-maximization harm the environment?},
  author={Wang, Lan and Cai, Gangshu and Tsay, Andy A and Vakharia, Asoo J},
  journal={Production and Operations Management},
  volume={26},
  number={8},
  pages={1585--1603},
  year={2017},
  publisher={SAGE Publications Sage CA: Los Angeles, CA}
}

@article{atasu2008remanufacturing,
  title={Remanufacturing as a marketing strategy},
  author={Atasu, Atalay and Sarvary, Miklos and Van Wassenhove, Luk N},
  journal={Management Science},
  volume={54},
  number={10},
  pages={1731--1746},
  year={2008},
  publisher={Informs}
}

@article{ferrer2006managing,
  title={Managing new and remanufactured products},
  author={Ferrer, Geraldo and Swaminathan, Jayashankar M},
  journal={Management Science},
  volume={52},
  number={1},
  pages={15--26},
  year={2006},
  publisher={INFORMS}
}

@article{gallego1994optimal,
  title={Optimal dynamic pricing of inventories with stochastic demand over finite horizons},
  author={Gallego, Guillermo and Van Ryzin, Garrett},
  journal={Management Science},
  volume={40},
  number={8},
  pages={999--1020},
  year={1994},
  publisher={INFORMS}
}

@article{li2024should,
  title={Should original equipment manufacturers authorize third-party remanufacturers?},
  author={Li, Wei and Jin, Mingzhou and Galbreth, Michael R},
  journal={European Journal of Operational Research},
  volume={314},
  number={3},
  pages={1013--1028},
  year={2024},
  publisher={Elsevier}
}

@article{banerjee2023optimal,
  title={Optimal patent licensing—Two or three-part tariff},
  author={Banerjee, Swapnendu and Mukherjee, Arijit and Poddar, Sougata},
  journal={Journal of Public Economic Theory},
  volume={25},
  number={3},
  pages={624--648},
  year={2023},
  publisher={Wiley Online Library}
}

@techreport{bonnet2006two,
  title = {Two-Part Tariffs versus Linear Pricing Between Manufacturers and Retailers: Empirical Tests on Differentiated Products Markets},
  author = {Christoph Carnehl and Anton Sobolev and Konrad Stahl and André Stenzel},
  institution = {CEPR Press},
  type = {CEPR Discussion Paper},
  number = {DP6016},
  year = {2006}}

@article{huang2024remanufacturing,
  title={Remanufacturing in global supply chains: Self-operating or licensing?},
  author={Huang, Hongfu and Xu, Fei and Wang, Min and Yang, Hui and Li, Taixin},
  journal={Transportation Research Part E: Logistics and Transportation Review},
  volume={190},
  pages={103708},
  year={2024},
  publisher={Elsevier}
}

@article{xiong2013don,
  title={Don’t forget your supplier when remanufacturing},
  author={Xiong, Yu and Zhou, Yu and Li, Gendao and Chan, Hingkai and Xiong, Zhongkai},
  journal={European Journal of Operational Research},
  volume={230},
  number={1},
  pages={15--25},
  year={2013},
  publisher={Elsevier}
}

@article{jin2023right,
  title={Right to repair: pricing, welfare, and environmental implications},
  author={Jin, Chen and Yang, Luyi and Zhu, Cungen},
  journal={Management Science},
  volume={69},
  number={2},
  pages={1017--1036},
  year={2023},
  publisher={INFORMS}
}

@article{atasu2010so,
  title={So what if remanufacturing cannibalizes my new product sales?},
  author={Atasu, Atalay and Guide Jr, V Daniel R and Van Wassenhove, Luk N},
  journal={California Management Review},
  volume={52},
  number={2},
  pages={56--76},
  year={2010},
  publisher={SAGE Publications Sage CA: Los Angeles, CA}
}

@techreport{WorldBank2022,
  title = {Squaring the Circle: Policies from Europe’s Circular Economy Transition},
  author = {{World Bank}},
  institution = {World Bank},
  year = {2022},
  page = {165}
}

@misc{units4.1,
title={Circular guide:
Remanufacturing Units 4.1},
author={{Circularity}},
year = {2022},
howpublished = {\url{https://circularity.com/en/circularguide/remanufacturing/#contactform}},
note={Accessed October 8, 2025}
}

@article{ginsburg2001once,
  title={Once is not enough},
  author={Ginsburg, Janet},
  journal={Business Week},
  number={3728},
  pages={128B--128D},
  year={2001},
  publisher={Proquest}
}

@book{chatti2019cirp,
  title={CIRP encyclopedia of production engineering},
  author={Chatti, Sami and Laperri{\`e}re, Luc and Reinhart, Gunther and Tolio, Tullio and others},
  year={2019},
  publisher={Springer}
}

@article{savaskan2004closed,
  title={Closed-loop supply chain models with product remanufacturing},
  author={Savaskan, R Canan and Bhattacharya, Shantanu and Van Wassenhove, Luk N},
  journal={Management Science},
  volume={50},
  number={2},
  pages={239--252},
  year={2004},
  publisher={INFORMS}
}

@article{du2012integrated,
  title={An integrated method for evaluating the remanufacturability of used machine tool},
  author={Du, Yanbin and Cao, Huajun and Liu, Fei and Li, Congbo and Chen, Xiang},
  journal={Journal of Cleaner Production},
  volume={20},
  number={1},
  pages={82--91},
  year={2012},
  publisher={Elsevier}
}

@article{san2015optimal,
  title={Optimal two-part tariff licensing mechanisms},
  author={San Mart{\'\i}n, Marta and Saracho, Ana I},
  journal={The Manchester School},
  volume={83},
  number={3},
  pages={288--306},
  year={2015},
  publisher={Wiley Online Library}
}

@article{arrow1951optimal,
  title={Optimal inventory policy},
  author={Arrow, Kenneth J and Harris, Theodore and Marschak, Jacob},
  journal={Econometrica: Journal of the Econometric Society},
  pages={250--272},
  year={1951},
  publisher={JSTOR}
}

@article{littlewood1972forecasting,
  title={Forecasting and control of passenger bookings},
  author={Littlewood, Kenneth},
  journal={The Airline Group of the International Federation of Operational Research Societies},
  volume={12},
  pages={95--117},
  year={1972}
}

@article{petruzzi1999pricing,
  title={Pricing and the newsvendor problem: A review with extensions},
  author={Petruzzi, Nicholas C and Dada, Maqbool},
  journal={Operations Research},
  volume={47},
  number={2},
  pages={183--194},
  year={1999},
  publisher={INFORMS}
}

\begin{appendix}
\numberwithin{equation}{section}
\setcounter{equation}{0} % Numbering from 0
\renewcommand{\theequation}{A\arabic{equation}}
\setcounter{table}{0}   % Numbering from 0
\setcounter{figure}{0}
\renewcommand{\thetable}{A\arabic{table}}
\renewcommand{\thefigure}{A\arabic{figure}}
\linespread{1.2}
\section*{Appendix}
We take Model N as a representative case for the following proofs; the same reasoning applies analogously to the other models.

\section{Compound Poisson demand}
Let the total number of potential consumers $N$ follows a Poisson distribution with rate $\lambda$, $\mathbb{P}(N=i)=e^{-\lambda} \frac{(\lambda)^i }{i!}.$ Each customer is characterized by a preference parameter $\theta \sim U[0,1]$ and will buy a new product if the utility condition $U_n(\theta)=\theta V_n-p_n \geq 0$ is met, that is, $\theta\geq\frac{p_n}{V_n}$. Therefore, the probability that a customer buys the product is $1-\frac{p_n}{V_n}$. 

Conditional on $N=i$, the number of customers who buy the new product follows the binomial distribution $\text{Binom}(i,1-\frac{p_n}{V_n})$.
By integrating over the distribution of $N$, demand $D_n$ follows a compound Poisson distribution:
\begin{align*}
\mathbb{P}(D_n=k) &= \sum_{i=k}^{\infty}\underbrace{e^{-\lambda}\frac{(\lambda)^i}{i!}}_{\text{Poisson probability}}\underbrace{\binom{i}{k}\left(1-\frac{p_n}{V_n}\right)^k\left(\frac{p_n}{V_n}\right)^{i-k}}_{\text{Binomial probability}} \notag \\
&= \frac{1}{k!}e^{-(1-\frac{p_n}{V_n})\lambda} \left[ (1-\frac{p_n}{V_n})\lambda\right]^k, \quad \forall k \in \mathbb{N}.
\label{eq:pdf demand n}
\end{align*}
The cumulative distribution function (CDF) for the demand $D_n$ is $F(k,\Lambda(p_n))=\sum_{i=0}^k \frac{1}{i!}e^{-\Lambda(p_n)} \Lambda^i(p_n)$. Hence, $D_n\sim \text{Poisson}\left((1-\frac{p_n}{V_n})\lambda\right)$. 

\section{Optimal quantity for any given price}
\label{app: model n optimal quantity}
Given the discrete nature of quantity $q_n$, we calculate the first and second differences to determine the optimal quantity for a given price. 

First, we derive the first difference $\Delta\Pi^N(p_n,q_n)$ for all $p_n$ and $q_n\in \mathbb{N}$. By the definition of the first difference, we have $\Delta\Pi^N(p_n,q_n) = \Pi^N(p_n,q_n+1)-\Pi^N(p_n,q_n).$

Substituting the expression for $\Pi^N(p_n,q_n)$ and simplifying yields 
\begin{align*}
\Delta\Pi^N(p_n,q_n) = p_n e^{-\Lambda(p_n)} \Bigg(\underbrace{\sum_{k=1}^{q_n+1}\frac{\Lambda^k(p_n)}{(k-1)!}-\sum_{k=1}^{q_n}\frac{\Lambda^k(p_n)}{(k-1)!}} + \underbrace{\sum_{k=q_n+2}^{\infty}(q_n+1) \frac{\Lambda^k(p_n)}{k!} - \sum_{k=q_n+1}^{\infty}q_n \frac{\Lambda^k(p_n)}{k!}}\Bigg) -c .
\end{align*}

Since the first pair of sums differ only by the term at $k=q_n+1$, we obtain $$\sum_{k=1}^{q_n+1}\frac{\Lambda^k(p_n)}{(k-1)!}-\sum_{k=1}^{q_n}\frac{\Lambda^k(p_n)}{(k-1)!}=\frac{\Lambda^{q_n+1}(p_n)}{q_n!}.$$

Isolating the overlapping range in the second pair of sums gives $$(q_n+1)\frac{\Lambda^{q_n+1}(p_n)}{(q_n+1)!} - q_n \frac{\Lambda^{q_n+1}(p_n)}{(q_n+1)!}=\frac{\Lambda^{q_n+1}(p_n)}{(q_n+1)!}.$$

Combining terms, we find $$\Delta\Pi^N(p_n,q_n) =p_n e^{-\Lambda(p_n)} \left(\frac{\Lambda^{q_n+1}(p_n)}{(q_n+1)!} + \sum_{k=q_n+2}^{\infty}\frac{\Lambda^k(p_n)}{k!}\right) - c.$$

Finally, merging the sums yields the desired results: 
\begin{equation*}
\Delta \Pi^N(p_n, q_n)
  = p_n e^{-\Lambda(p_n)}\sum_{k=q_n+1}^{\infty}\frac{\Lambda^k(p_n)}{k!} - c.
\end{equation*}

It is straightforward to verify that $$\Delta^2 \Pi^N(p_n,q_n)= -p_n e^{-\Lambda(p_n)}\frac{\Lambda^{q_n+1}(p_n)}{(q_n+1)!}.$$
Given that $p_n>0$ and $q_n\in\mathbb{N}$, it immediately follows that $\Delta^2 \Pi^N(p_n,q_n)<0$ for all feasible $q_n$. Hence, the profit function $\Pi^N(p_n,q_n)$ is strictly concave with respect to $q_n$. 

Using the summation in terms of the CDF of the Poisson demand $F(k,\Lambda(p_n))$, the first difference can be equivalently written as $$\Delta\Pi^N(p_n,q_n)=e^{-\Lambda(p_n)}\left[1-F(q_n,\Lambda(p_n)\right]-c.$$ 

Setting $\Delta\Pi^N(p_n,q_n)=0$ yields the optimal quantity for any given price $q_n^*(p_n)=F^{-1}(1-\frac{c}{p_n}).$

\section{Reduced-form optimization problem}
Recall that the OEM's expected profit function from the new product, given price $p_n$ and quantity $q_n$, is defined as $$\Pi^N(q_n,p_n)= p_n e^{-\Lambda(p_n)} \left[\sum_{k=1}^{q_n}\frac{\Lambda^k(p_n)}{(k-1)!} + \sum_{k=q_n+1}^{\infty} q_n \frac{\Lambda^k(p_n)}{k!} \right] -c q_n.$$

Substituting the optimal quantity $q_n^*(p_n)=F^{-1}(1-\frac{c}{p_n})$ yields 
\begin{align*}
\Pi^N(p_n,q_n^*)
&=p_n e^{-\Lambda(p_n)} \left[ \sum_{k=1}^{q_n^*}  \frac{\Lambda^k(p_n)}{(k-1)!} + q_n^*\sum_{k=q_n^*+1}^{\infty}\frac{\Lambda^k(p_n)}{k!}\right] - c q_n^*\\
&=p_n \left[ \sum_{k=1}^{q_n^*} e^{-\Lambda(p_n)} \frac{\Lambda^k(p_n)}{(k-1)!} +q_n^*\text{ Pr}(Dn>q^*_n)\right] - c q_n^*.
\end{align*}
where $\text{Pr}(D_n=k)=e^{-\Lambda(p_n)} \frac{\Lambda^k(p_n)}{(k)!}$ and $\text{Pr}(D_n>q^*_n)=1-F(q^*_n,\Lambda(p_n))$.

For the first summation, we have $\sum_{k=1}^{q_n^*} e^{-\Lambda(p_n)} \frac{\Lambda^k(p_n)}{(k-1)!}=\Lambda(p_n)\sum_{k=1}^{q_n^*} e^{-\Lambda(p_n)} \frac{\Lambda^(k-1)(p_n)}{(k-1)!}= \Lambda(p_n) F(q^*_n-1,\Lambda(p_n))$ since the inner sum is the CDF of a Poisson variable evaluated at $q_n^*-1$.

Substituting and using the optimal condition $1- F(q_n^*,\Lambda(p_n))=\frac{c}{p_n}$ gives $$\Pi^N(p_n,q_n^*)=p_n \Lambda(p_n)F(q_n^* -1,\Lambda(p_n)).$$ 

Thus, the optimization problem reduces to $\max_{p_n} \ \Pi^N(p_n,q_n^*)
= p_n \Lambda(p_n) F(q_n^* -1,\Lambda(p_n)).$

\section{Jointly demand rates $\texorpdfstring{\Lambda_n}{Lambdan}$ and $\texorpdfstring{\Lambda_r}{Lambdar}$}
\label{app: Model O poisson demands}
Similarly to Model N, we derive the purchasing probability for new and remanufactured products by analyzing the customer utilities. 

A consumer with preference $\theta$ obtains utility $U_n(\theta) = \theta \left(\delta V + \beta(\delta V-\alpha \delta V)\right) - p_{n}$, $U_r(\theta) = \theta \alpha \delta V - p_r$, and $U_0(\theta)=0$ for buying a new product, a remanufactured products and nothing, respectively. 

A customer prefers to buy a new product if $U_n(\theta)\geq U_r(\theta)$ and $U_n(\theta)\geq 0$. Solving these conditions gives $\theta \geq \max\{\frac{p_{n}-p_r}{(1+\beta)(1-\alpha)\delta V }, \frac{p_{n}}{ (1+\beta -\alpha \beta)\delta V}\}$.
Therefore, the probability of buying a new product is 
\begin{align}
\text{Prob}_{n}(p_n,p_r) &:= 1-\min\left\{1,\max\left\{\frac{p_{n}-p_r}{(1+\beta)(1-\alpha)\delta V}, \frac{p_{n}}{ (1+\beta -\alpha \beta)\delta V}\right\}\right\}, \notag\\
&=\begin{cases}
1-\frac{p_{n}-p_r}{(1+\beta)(1-\alpha)\delta V}, & \text{ if } \frac{p_{n}}{(1+\beta -\alpha \beta)\delta V} \leq \frac{p_{n}-p_r}{(1+\beta)(1-\alpha)\delta V} \leq 1;\\
1-\frac{p_{n}}{(1+\beta -\alpha \beta)\delta V}, & \text { if } \frac{p_{n}-p_r}{(1+\beta)(1-\alpha)\delta V } \leq \frac{p_{n}}{(1+\beta -\alpha \beta)\delta V} \leq 1;\\
0, & \text { otherwise }.
\end{cases}
\label{eq: Model O prob n}
\end{align}

A consumer prefers the remanufactured product if $U_r(\theta) > U_n(\theta)$ and $U_r(\theta)\geq 0$, leading to $\theta \geq \frac{p_{r}}{\alpha \delta V}$ and $\theta <\frac{p_{n}-p_r}{(1+\beta)(1-\alpha)\delta V}$. Hence, the probability of buying a remanufactured product is: 
\begin{align}
\text{Prob}_{r}(p_n,p_r) &:= \max\left\{0, \min\left\{1,\frac{p_{n}-p_r}{(1+\beta)(1-\alpha)\delta V}\} - \min\{1,\frac{p_{r}}{\alpha \delta V}\right\}\right\}, \notag \\
&=\begin{cases}
\frac{p_{n}-p_r}{(1+\beta)(1-\alpha)\delta V} - \frac{p_{r}}{\alpha \delta V}, & \text{ if } \frac{p_{r}}{\alpha \delta V}\leq \frac{p_{n}-p_r}{(1+\beta)(1-\alpha)\delta V} \leq 1;\\
1- \frac{p_{r}}{\alpha \delta V}, & \text { if } \frac{p_{r}}{\alpha \delta V}\leq 1 \leq \frac{p_{n}-p_r}{(1+\beta)(1-\alpha)\delta V};\\
0, & \text { otherwise }.
\end{cases}
\label{eq: Model O prob r}
\end{align}

According to the probabilities given in Eqs.~\eqref{eq: Model O prob n} and ~\eqref{eq: Model O prob r}, the demands for new and remanufactured products, $D_n$ and $D_r$, follow Poisson distributions with demand rates $\Lambda_n(p_n,p_r) : = \text{Prob}_{n}(p_n,p_r) \lambda$ and $\Lambda_r(p_n,p_r) : = \text{Prob}_{r}(p_n,p_r) \lambda$, respectively. 

\section{Demand region charts}
\label{app: regions}
To examine how the price pair $(p_{n},p_r)$ influences market demands, we analyze the corresponding demand rates $\Lambda_{n}$ and $\Lambda_r$ under different conditions. 
We define the following threshold parameters derived from Eqs.~\eqref{eq: Model O prob n} and \eqref{eq: Model O prob r}:
\begin{equation*}
    \hat{\theta}_1 := \frac{p_{n}-p_r}{(1+\beta)(1-\alpha)\delta V}, \quad \hat{\theta}_2 := \frac{p_{n}}{(1+\beta -\alpha \beta)\delta V},
    \quad \hat{\theta}_3 := \frac{p_{r}}{\alpha \delta V}.
\end{equation*}

By combining the conditions implied by $\hat{\theta}_1$, $\hat{\theta}_2$ and $\hat{\theta}_3$, we summarize all cases of the value of demand rates in Table~\ref{table: Model O joint conditions}.
We reformulate these conditions of $\hat{\theta}_1$, $\hat{\theta}_2$, and $\hat{\theta}_3$ as relationships between $p_n$ and $p_r$, and display identify distinct market regions in Figs.~\ref{fig: Model O demand rates cases} and ~\ref{fig: Model T regions}.

\begin{table}[ht!]\small
\centering
\caption{Joint conditions for $\Lambda_n(p_n,p_r)$ and $\Lambda_r(p_n,p_r)$} \label{table: Model O joint conditions}
\begin{tabular}{|c|c|c|c|c|}
    \hline 
    \textbf{Cases} & \multicolumn{2}{c|}{\textbf{Conditions}} & \textbf{$\Lambda_n(p_n,p_r)$} & \textbf{$\Lambda_r(p_n,p_r)$} \\
    \hline 
    (1) &
    \multirow{2}{*}{$\hat{\theta}_1 < \hat{\theta}_2 < \hat{\theta}_3$} 
    & $\hat{\theta}_2>1$ 
    & 0 
    & 0 \\
    (2) &
    & $\hat{\theta}_2 < 1$ 
    & $\left(1 - \frac{p_{n}}{(1+\beta -\alpha \beta)\delta V}\right)  \lambda$ 
    & 0 \\
    \hline 
    (3) &
    \multirow{3}{*}{$\hat{\theta}_3 < \hat{\theta}_2 < \hat{\theta}_1$} 
    & $1 < \hat{\theta}_3$ 
    & 0 
    & 0 \\
    (4) &
    & $\hat{\theta}_3< 1 < \hat{\theta}_1$ 
    & 0 
    & $\left(1-\frac{p_{r}}{\alpha \delta V}\right) \lambda$ \\
    (5) &
    & $\hat{\theta}_1<1$
    & $\left(1-\frac{p_{n}-p_r}{(1+\beta)(1-\alpha)\delta V}\right) \lambda$ 
    & $\left(\frac{p_{n}-p_r}{(1+\beta)(1-\alpha)\delta V} - \frac{p_{r}}{\alpha \delta V}\right) \lambda$ \\
    \hline
\end{tabular}
\end{table}

\section{Approximate boundary \texorpdfstring{$\beta_1(\alpha)$}{beta1(alpha)}}
Replacing $F^*(q^*-1)$ with $F(q^*)$ in Eqs.~\eqref{eq: model n expected profit with optimal q} and \eqref{eq: Model O profit with optimal q}, we obtain the approximate profit functions $\widetilde{\Pi}^N(p_n) = (p_n-c) \Lambda(p_n),$ and $\widetilde{\Pi}^O(p_{n},p_r) = (p_{n}-c) \Lambda_{n}+(p_{r}-c_r-c_{coll}) \Lambda_{r}.$ These follow the classical Newsvendor formulation, where $\Lambda$ represents expected sales. Solving the first- and second-order conditions, we yield closed-form approximation solutions. 

For Model N, the approximate price and profit are $\widetilde{p_n}^{N*} = \frac{c+\delta V}{2}$ and $\ \widetilde{\Pi}^{N*}(\cdot) = \frac{(c-\delta V)^2}{4 \delta V}\lambda$, respectively.

For Model O, three market equilibria emerge. 
\begin{itemize}
    \item Coexistence market: When $\frac{(1+\beta^-)(1-\alpha)}{\alpha}(c_r+c_{coll})<c-c_r-c_{coll}< (1+\beta^-)(1-\alpha)\delta V$, the approximate prices are $\widetilde{p_n}^{O*} = \frac{c+(1+\beta^--\alpha \beta^-)\delta V}{2},\widetilde{p_r}^{O*} =\frac{c_r+c_{coll}+\alpha \delta V}{2}$. The corresponding approximate expected profit is $
    \widetilde{\Pi}^{O*} = 
    \frac{\lambda}{4 \delta V} \Big[
    \frac{\big((1+\beta^-) (1-\alpha) \delta V-(c+c_r+c_{coll})\big)^2}{(1+\beta^-) (1-\alpha)} + \frac{\big( \alpha\delta V - c_r-c_{coll}\big)^2}{\alpha} \Big].$
    \item Only remanufactured market: When $c-c_r-c_{coll}\geq (1+\beta^-)(1-\alpha)\delta V$, the approximate price is $\widetilde{p_r}^{O*} = \frac{c_r+c_{coll}+\alpha\delta V}{2}$, and the corresponding approximate expected profit is $\widetilde{\Pi}^{O*} = \frac{\lambda}{4 \alpha \delta V}  \big(c_r+c_{coll}-\alpha \delta V\big)^2$.
    \item Only new market: When $c-c_r-c_{coll}\leq \frac{(1+\beta^-)(1+\alpha)}{\alpha}(c_r+c_{coll})$, the equilibrium market is dominated by new products, which is equivalent to Model N.
\end{itemize}

By comparing the approximate profits of Models N and O, we solve the closed-form approximate boundary is derived as:{\small$$\beta_1(\alpha)=|\frac{\alpha c^2-(c_r+c_{coll})^2-\alpha(1-\alpha)(\delta V)^2+\sqrt{\left( (\alpha \delta V+c_r+c_{coll})^2-\alpha(\delta V+c)^2 \right) \cdot \left( (\alpha \delta V-c_r-c_{coll})^2-\alpha(\delta V-c)^2 \right)}}{2\alpha(1-\alpha)(\delta V)^2}|.$$}
The function is defined over $[\alpha_1',\alpha_2]$, where $\alpha_1' = \frac{c_r+c_{coll}}{c}$, $\alpha_2= \frac{(\delta V-c)^2 + (\delta V-c)\sqrt{(\delta V-c)^2+4\delta V (c_r+c_{coll})}+2\delta V (c_r+c_{coll})}{2 \delta^2 V^2}$.

\end{appendix}

\end{document}

\begin{center}
    {\LARGE Supplementary Material}
\end{center}

%%%%%%%%%%%%%%%%%%%%%%%%%%%%%%%%
% Analytical proofs
\section{Proofs} \label{app: proof}

\subsection{Model N: Demand analysis}
\label{app: model n poisson demand}

% \begin{lemma}
%     Without remanufacturing, the demand for new products follows a Poisson distribution $D_n \sim \text{Poisson}(\Lambda(p_n))$, with the price-sensitive demand rate $\Lambda(p_n) := (1-\frac{p_n}{V_n}) \cdot \lambda$. The cumulative distribution function (CDF) is given by $F(k,\Lambda(p_n)):=\sum_{i=0}^{k}e^{-\Lambda(p_n)}\frac{\Lambda^{-i}(p_n)}{i!}, \quad \forall k \in \mathbb{N}$.
%     \label{lemma: model n poisson demand}
% \end{lemma}
The customer utility function determines whether a customer with preference parameter $\theta_i$, $i \in \{0,1,\cdots,N\}$ will make a purchase based on 
\begin{align*}
\mathbbm{1}_{U_n(\theta_i)\geq 0}(\theta_i):= 
\begin{cases}
1, & \text { if } U_n(\theta_i)=\theta_i V_n-p_n \geq 0;\\
0, & \text { if } U_n(\theta_i)=\theta_i V_n-p_n <0,
\end{cases}
\end{align*} where $\mathbbm{1}_{U_n(\theta_i)\geq 0}(\theta_i)=1$ represents a purchase decision and $0$ represents no purchase. Since customers buy when $\theta_i V_n-p_n \geq 0$, we can determine that purchases occur when  $\theta_i\geq\frac{p_n}{V_n}$. Given that $\theta$ follows a uniform distribution on $[0,1]$, each customer has a probability of $(1-\frac{p_n}{V_n})$ of making a purchase. 

We assume that the total number of potential consumers $N$ follows a Poisson distribution with parameter $\lambda$, that is $\mathbb{P}(N=i)=e^{-\lambda} \frac{(\lambda)^i }{i!}$. For any given number of arrivals $i$, each customer independently decides to buy a new product with probability $(1-\frac{p_n}{V_n})$, which generates a binomial distribution with probability parameters $n$ and $(1-\frac{p_n}{V_n})$. 

Combining those two parts, the demand of new products  $D_n$ follows a compound Poisson distribution, where we sum over all possible numbers of arrivals: 
\begin{align}
\mathbb{P}(D_n=k) &= \sum_{i=k}^{\infty}\underbrace{e^{-\lambda}\frac{(\lambda)^i}{i!}}_{\text{Poisson probability}}\underbrace{\binom{i}{k}\left(1-\frac{p_n}{V_n}\right)^k\left(\frac{p_n}{V_n}\right)^{i-k}}_{\text{Binomial probability}} \notag \\
&= \frac{1}{k!}e^{-(1-\frac{p_n}{V_n})\lambda} \left[ (1-\frac{p_n}{V_n})\lambda\right]^k, \quad \forall k \in \mathbb{N}.
\label{eq:pdf demand n}
\end{align}
This expression shows that the demand $D_n$ follows a Poisson distribution with parameter $(1-\frac{p_n}{V_n})\lambda$. For simplicity, we construct an auxiliary variable as $\Lambda(p_n):=(1-\frac{p_n}{V_n})\lambda$, which denotes the Poisson demand rate. Then the cumulative distribution function can be calculated as $F(k,\Lambda(p_n)):=\mathbb{P}(D_n\leq k)=\sum_{i=0}^k \mathbb{P}(D_n=i) =\sum_{i=0}^k \frac{1}{i!}e^{-\Lambda(p_n)} \Lambda^i(p_n), \quad \forall k \in \mathbb{N}.$

\subsection{Model N: Approximation method}
\label{app: model n same solution}
First we complete the derivation of approximate function. 
Recall that the expected profit function is given by $\Pi^N(p_n,q^*_n(p_n)) = p_n \cdot \Lambda(p_n) \cdot F\left(q^*(p_n)-1, \Lambda(p_n)\right)$, where $q_n^*(p_n)=F^{-1}(1-\frac{c}{p_n})$. By substituting the demand rate $\Lambda(p_n)=(1-\frac{p_n}{V_n})\lambda$, we obtain: 
\begin{equation}
    \Pi^N(p_n,q^*_n(p_n)) = p_n \cdot (1-\frac{p_n}{V_n}) \lambda \cdot F\left(q_n^*(p_n)-1, \Lambda(p_n)\right).
    \label{eq: model n profit with F(q-1)}
\end{equation}

Our approximation approach replaces $F\left(q^*_n(p_n)-1,\Lambda(p_n)\right)$ with $F\left(q^*_n(p_n),\Lambda(p_n)\right)$, therefore the approximate expected profit can be expressed as:
\begin{align}
\widetilde{\Pi}^N(p_n,q^*_n(p_n)) 
&= p_n \cdot \left(1-\frac{p_n}{V_n}\right) \lambda \cdot F\left(q^*_n(p_n),\Lambda(p_n)\right) \notag \\
& =  p_n \cdot \left(1-\frac{p_n}{V_n}\right) \lambda \cdot \left(1-\frac{c}{p_n}\right) \notag \\
& = \left(p_n-c\right) \cdot \left(1-\frac{p_n}{V_n}\right)\lambda. 
\label{eq: model n approximate function}
\end{align}

To prove $\widetilde{\Pi}^N(p_n,q^*_n(p_n))$ connects all local maximum points of the original function $\Pi^N(p_n,q^*_n(p_n))$, we first define the set of local maximum points of the original function, and then prove within this set, the approximate function yields the same optimal solution as the original function. 

Due to the discrete properties of ICDF $q_n^*(p_n)=F^{-1}(1-\frac{c}{p_n})$, there exists intervals where $q_n^*(p_n)$ maintains a constant value within each interval. To simplify the formulation, we define an interval $[p_n^-, p_n^+)$ as an example to explain the property of the original profit function. Within this interval, $F(q_n^*(p_n), \Lambda(p_n))$ and $F(q_n^*(p_n)-1, \Lambda(p_n))$ are monotonically increasing with $p_n$. This results in $\Pi^N(p_n, q_n^*(p_n))$ monotonically increasing within this interval. When $p_n$ is less than $p_n^-$, $q_n^*(p_n)$ increases by one unit; when $p_n$ exceeds $p_n^+$, $q_n^*(p_n)$ decreases by one unit, creating a sharp downwards and discontinuity. Thus, on a larger scale, the monotonic increases within intervals and sudden drops at the boundary of intervals create a sawtooth pattern. Note that the left limit of $p_n^+$  behaves as a local maximum point.

Based on the above analysis, we define the support set for the discrete optimal quantities: $\mathcal{Q} := \{q^k(p_n): q^k = F^{-1}(1-\frac{c}{p_n}), q^k>q^{k+1}, p_n \in [c, \delta V], k \in \mathbb{N} \}$, with $q^0 = F^{-1}(1-\frac{c}{\underline{p}_n})$. Additionally, we construct two sets for prices based on discontinuities: 
\begin{enumerate}
    \item Starting Points: $\mathcal{SP} = \{p_n^{s,k}: p_n^{s,k} := \inf\{p_n: q_n^*(p_n)=q^k, q^k \in \mathcal{Q}\}\}$,
    \item Ending Points: $\mathcal{EP} = \{p_n^{e,k}:p_n^{e,k} := \sup\{p_n: q_n^*(p_n)=q^k, q^k \in \mathcal{Q}\}\}$.
\end{enumerate}

That is, $p_n^{s,k}$ represents the greatest lower bound (included), and $p_n^{e,k}$ represents the least upper bound (not included) of the interval that keeps $q_n^*(p_n)$ remain $q^k$. This can be expressed mathematically as $q^k = q^*(p_n^{s,k}) = \lim_{\epsilon \rightarrow 0^+}q^*(p_n^{e,k}-\epsilon)$. One may understand that for any given $k$, the interval $[p_n^{s,k},p_n^{e,k})$ corresponds to the generic interval $[p_n^+,p_n^-)$ of the simple example. And the ending points in set $\mathcal{EP}$ are the local maximum points. The left-closed, right-open structure of these intervals creates a crucial relationship: The least upper bound of interval $k$ ($p_n^{e,k}$) equals the greatest lower bound of interval $k+1$ ($p_n^{s,k+1}$). This can be expressed mathematically as $p_n^{s,k+1} = p_n^{e,k}$.
Because of the continuity of interval, optimal quantities shall change one unit at a time, that is, $q^k = q^{k+1}+1$.

In conclusion, the following relations hold:
\begin{align*}
&p_n^{s,k+1} = p_n^{e,k}; \\
&q^k = q^*(p_n^{s,k}) = \lim_{\epsilon \rightarrow 0^+}q^*(p_n^{e,k}-\epsilon);\\
&q^k = q^{k+1}+1.
\end{align*}
Figure~\ref{fig:price sets} illustrates an example of the original profit function and optimal quantity with respect to the sets of optimal quantities and prices.
\begin{figure}
\centering
\includegraphics[width=0.7\linewidth]{images/appendix_A5.pdf}
\caption{Illustration of sets $\mathcal{Q}$ and $\mathcal{SP},\mathcal{EP}$} \label{fig:price sets}
\end{figure}

From the definition of $\mathcal{Q}$, we have $F(q^k,\Lambda(p_n^{s,k})) = 1-\frac{c}{p_n^{s,k}}, k\in \mathbb{N}$, and derive that $F(q^{k-1}-1,\Lambda(p_n^{e,k-1})) = F(q^k,\Lambda(p_n^{s,k})) = 1-\frac{c}{p_n^{s,k}}, k\in \mathbb{N^+}$. Using this expression, we reformulate the original profit function, defined by Eq.~\eqref{eq: model n profit with F(q-1)}, on the support set $\mathcal{EP}$ as follows:
\begin{align*}
    \max_{p_n \in \mathcal{EP}} \Pi^N\left(p_n, q^*_n(p_n)\right) &= \Pi^N(p_n^{e,k},q^k) \notag \\
    & = p_n^{e,k} \cdot \left(1-\frac{p_n^{e,k}}{V_n}\right) \lambda \cdot F\left(q^k-1, \Lambda(p_n^{e,k})\right) \notag \\
    &= p_n^{e,k} \cdot \left(1-\frac{p_n^{e,k}}{V_n}\right) \lambda \cdot \left(1-\frac{c}{p_n^{s,k+1}}\right) \notag \\ 
    &= p_n^{e,k} \cdot \left(1-\frac{p_n^{e,k}}{V_n}\right) \lambda \cdot \left(1- \frac{c}{p_n^{e,k}}\right) \notag \\ 
    &= (p_n^{e,k}-c) \cdot \left(1-\frac{p_n^{e,k}}{V_n}\right) \lambda.
    % \label{eq: model n subproblem with EP}
\end{align*}

Comparing it with the approximate profit function from Eq.~\eqref{eq: model n approximate function}, we prove that $\Pi^N(\cdot)$ and $\widetilde{\Pi}^N(\cdot)$ coincide at the points in $\mathcal{EP}$, which behave as the local maximum points of the original function.

%%%%%%%%%%%%%%%%%%%%%%%%%%%%%%%%
% Model O: In-house Remanufacture
%%%%%%%%%%%%%%%%%%%%%%%%%%%%%%%%
% \section{Model O: In-house remanufacture}
% \label{app: Model O}
% \setcounter{equation}{0} % Numbering from 0
% \renewcommand{\theequation}{B\arabic{equation}}
% \setcounter{table}{0}   % Numbering from 0
% \setcounter{figure}{0}
% \renewcommand{\thetable}{B\arabic{table}}
% \renewcommand{\thefigure}{B\arabic{figure}}

%%%%%%%%%%%%%% Poisson Demands %%%%%%%%%%%%%%%%%%
\subsection{Model O: Demand analysis}
\label{app: Model O poisson demands}
Similarly to Model N, we analyze the customer utilities to derive the possible density function (PDF) of demands $D_n$ and $D_r$. 

Recall that consumer's utility from buying a new product is $U_n(\theta) = \theta \left(\delta V + \beta(\delta V-\alpha \delta V)\right) - p_{n}$, that from buying a remanufactured product is $U_r(\theta) = \theta \alpha \delta V - p_r$, and that from buying nothing is $0$. 
Calculating $U_n(\theta)\geq U_r(\theta)$ and $U_n(\theta)\geq 0$, we deduce that customers with $\theta \geq \max\{\frac{p_{n}-p_r}{(1+\beta)(1-\alpha)\delta V }, \frac{p_{n}}{ (1+\beta -\alpha \beta)\delta V}\}$ will buy a new product to obtain higher surplus. Given that $\theta \sim \text{U}[0,1]$, each customer buys a new product with probability $\text{Prob}_{n}(p_n,p_r)$, where
\begin{align}
\text{Prob}_{n}(p_n,p_r) &:= 1-\min\{1,\max\{\frac{p_{n}-p_r}{(1+\beta)(1-\alpha)\delta V}, \frac{p_{n}}{ (1+\beta -\alpha \beta)\delta V}\}\}, \notag\\
&=\begin{cases}
1-\frac{p_{n}-p_r}{(1+\beta)(1-\alpha)\delta V}, & \text{ if } \frac{p_{n}}{(1+\beta -\alpha \beta)\delta V} \leq \frac{p_{n}-p_r}{(1+\beta)(1-\alpha)\delta V} \leq 1;\\
1-\frac{p_{n}}{(1+\beta -\alpha \beta)\delta V}, & \text { if } \frac{p_{n}-p_r}{(1+\beta)(1-\alpha)\delta V } \leq \frac{p_{n}}{(1+\beta -\alpha \beta)\delta V} \leq 1;\\
0, & \text { otherwise }.
\end{cases}
\label{eq: Model O prob n}
\end{align}

Given that the demand for new products is $D_n=\sum_{i=1}^{N}\mathbbm{1}_{(U_{n}\geq U_r, U_{n}\geq 0)}(\theta_i)$ and total customer of potential consumers $N$ follow Poisson distribution, we derive the PDF of demands $D_{n}$ as follows:
\begin{align}
\mathbb{P}(D_{n}=k) 
&=\sum_{i=k}^{\infty}e^{-\lambda} \cdot \frac{(\lambda)^i }{i!} \cdot \binom{i}{k} \cdot \text{Prob}_{n}^k \cdot (1-\text{Prob}_{n})^{i-k} \notag \\
&= \frac{1}{k!} \cdot e^{- \text{Prob}_{n} \lambda} \cdot \left(\text{Prob}_{n} \lambda\right)^k, \quad \forall k \in \mathbb{N}.
\label{eq: Model O pdf of demand n}
\end{align}

By calculating $U_r(\theta) > U_n(\theta)$ and $U_r(\theta)\geq 0$, we have that customers with $\theta \geq \frac{p_{r}}{\alpha \delta V}$ and $\theta <\frac{p_{n}-p_r}{(1+\beta)(1-\alpha)\delta V}$ will buy remanufactured products, and the probability of each customer will buy remanufactured product can be expressed as follows: 
\begin{align}
\text{Prob}_{r}(p_n,p_r) &:= \max\{0, \min\{1,\frac{p_{n}-p_r}{(1+\beta)(1-\alpha)\delta V}\} - \min\{1,\frac{p_{r}}{\alpha \delta V}\}\}, \notag \\
&=\begin{cases}
\frac{p_{n}-p_r}{(1+\beta)(1-\alpha)\delta V} - \frac{p_{r}}{\alpha \delta V}, & \text{ if } \frac{p_{r}}{\alpha \delta V}\leq \frac{p_{n}-p_r}{(1+\beta)(1-\alpha)\delta V} \leq 1;\\
1- \frac{p_{r}}{\alpha \delta V}, & \text { if } \frac{p_{r}}{\alpha \delta V}\leq 1 \leq \frac{p_{n}-p_r}{(1+\beta)(1-\alpha)\delta V};\\
0, & \text { otherwise }.
\end{cases}
\label{eq: Model O prob r}
\end{align}
The corresponding PDF of the demand for remanufactured products, denoted by $D_r=\sum_{i=1}^{N}\mathbbm{1}_{(U_r>U_{n}, U_r\geq 0)}(\theta_i)$, can be expressed as:
\begin{align}
\mathbb{P}(D_r=k) 
&=\sum_{i=k}^{\infty}e^{-\lambda} \cdot \frac{(\lambda)^i }{i!} \cdot \binom{i}{k} \cdot \text{Prob}_{r}^k \cdot (1-\text{Prob}_{r})^{i-k} \notag \\
&= \frac{1}{k!} \cdot e^{- \text{Prob}_{r} \lambda} \cdot \left(\text{Prob}_{r} \lambda\right)^k,\quad \forall k \in \mathbb{N}.
\label{eq: Model O pdf of demand r}
\end{align}

According to the PDF expressions in Eq.s~\eqref{eq: Model O pdf of demand n} and \eqref{eq: Model O pdf of demand r}, the demands $D_n$ and $D_r$ follows Poisson distributions with demand rates $\Lambda_n(p_n,p_r) : = \text{Prob}_{n}(p_n,p_r) \lambda$ and $\Lambda_r(p_n,p_r) : = \text{Prob}_{r}(p_n,p_r) \lambda$, respectively. For simplicity, we abbreviate them as $\Lambda_n$ and $\Lambda_r$, respectively in the following formulations. Note that the definition of demand rates occurs when both products are present on the market, which makes a comparison between $U_n$ and $U_r$ feasible; if only new or remanufactured products exist on the market, the demand rate should be derived as in single-product market.

\subsection{Derivation of Eq. (7)}
% ~\eqref{eq: Model O profit}}
\label{app:model o optimization problem}
Substituting the PDFs of demands into the expected sales $\mathbb{E}\left[\min\{D_{n},q_{n}\}\right], \mathbb{E}\left[\min\{D_{r},q_{r}\}\right]$, we obtain: 
\begin{alignat}{2}
    \mathbb{E}\left[\min\{D_{n},q_{n}\}\right] &= \sum_{k=0}^{q_{n}} k \cdot \mathbb{P}(D_{n}=k) + \sum_{k=q_{n}+1}^{\infty} q_{n} \cdot \mathbb{P}(D_{n}=k) \notag\\
    &= \sum_{k=0}^{q_{n}} k \cdot \frac{1}{k!} \cdot e^{-\text{Prob}_{n} \lambda} \cdot (\text{Prob}_{n} \lambda)^k + \sum_{k=q_{n}+1}^{\infty} q_{n} \cdot \frac{1}{k!} \cdot e^{-\text{Prob}_{n} \lambda} \cdot (\text{Prob}_{n} \lambda)^k \notag \\
    &= e^{-\text{Prob}_{n} \lambda} \cdot \left[\sum_{k=1}^{q_{n}} \frac{1}{(k-1)!} \cdot (\text{Prob}_{n} \lambda)^k + \sum_{k=q_{n}+1}^{\infty} q_{n} \cdot \frac{1}{k!} \cdot (\text{Prob}_{n} \lambda)^k\right] \notag \\
    &= e^{-\Lambda_{n}} \cdot \left[\sum_{k=1}^{q_{n}} \frac{1}{(k-1)!} \cdot \Lambda_{n}^k + \sum_{k=q_{n}+1}^{\infty} q_{n} \cdot \frac{1}{k!} \cdot \Lambda_{n}^k\right], \label{eq: Model O expected sales n}\\
    \mathbb{E}\left[\min\{D_r,q_r\}\right] &= \sum_{k=0}^{q_r} k \cdot \mathbb{P}(D_{r}=k) + \sum_{k=q_{r}+1}^{\infty} q_{r} \cdot \mathbb{P}(D_{r}=k) \notag \\
    &= e^{-\Lambda_{r}} \left[\sum_{k=1}^{q_{r}} \frac{1}{(k-1)!} \cdot \Lambda_{r}^k + \sum_{k=q_{r}+1}^{\infty} q_{r} \cdot \frac{1}{k!} \cdot \Lambda_{r}^k\right]. \label{eq: Model O expected sales r}
\end{alignat}

Recall that $\Pi^{\text{O}}(p_{n},q_{n},p_r,q_r) = p_{n}\cdot \mathbb{E}\big[\min \{D_n,q_n\}\big] + p_r \cdot \mathbb{E}\big[\min \{D_r,q_r\}\big] - c\cdot q_{n}- (c_r +c_{coll}) \cdot q_r$. By substituting Eq.~\eqref{eq: Model O expected sales n}-\eqref{eq: Model O expected sales r} into the expected profit function, we can obtain the optimization problem:
\begin{align*}
    \max_{p_{n},q_{n},p_r,q_r} \ \Pi^{\text{O}}(p_{n},q_{n},p_r,q_r) = 
    & \ p_{n}\cdot e^{-\Lambda_{n}} \cdot \left[\sum_{k=1}^{q_{n}} \frac{1}{(k-1)!} \cdot \Lambda_{n}^k+ \sum_{k=q_{n}+1}^{\infty} q_{n} \cdot \frac{1}{k!} \cdot \Lambda_{n}^k\right] \\
    &+ p_r \cdot e^{-\Lambda_{r}} \left[\sum_{k=1}^{q_{r}} \frac{1}{(k-1)!} \cdot \Lambda_{r}^k + \sum_{k=q_{r}+1}^{\infty} q_{r} \cdot \frac{1}{k!} \cdot \Lambda_{r}^k\right]\\
    & -c\cdot q_{n}-(c_r+c_{coll})\cdot q_r.
\end{align*}

\subsection{Explanation for Fig. 3}%~\ref{fig: Model O demand rates cases}}
\label{app: Model O regions}
To further investigate how $p_{n}$ and $p_r$ affect demands, we analyze the values of $\Lambda_{n}$ and $\Lambda_r$ and their conditions. We define the following terms to simplify the notations in Eq.s~\eqref{eq: Model O prob n} and \eqref{eq: Model O prob r}:
\begin{equation*}
    \hat{\theta}_1 := \frac{p_{n}-p_r}{(1+\beta)(1-\alpha)\delta V}, \quad \hat{\theta}_2 := \frac{p_{n}}{(1+\beta -\alpha \beta)\delta V},
    \quad \hat{\theta}_3 := \frac{p_{r}}{\alpha \delta V}.
\end{equation*}

By combining these conditional boundaries $\hat{\theta}_1$, $\hat{\theta}_2$ and $\hat{\theta}_3$, we derive the joint conditions for different values of $\Lambda_n(p_n,p_r)$ and $\Lambda_r(p_n,p_r)$ in Table~\ref{table: Model O joint conditions}.
\begin{table}[b]\small
\centering
\caption{Joint conditions for $\Lambda_n(p_n,p_r)$ and $\Lambda_r(p_n,p_r)$} \label{table: Model O joint conditions}
\begin{tabular}{|c|c|c|c|c|}
    \hline 
    \textbf{Cases} & \multicolumn{2}{c|}{\textbf{Conditions}} & \textbf{$\Lambda_n(p_n,p_r)$} & \textbf{$\Lambda_r(p_n,p_r)$} \\
    \hline 
    (1) &
    \multirow{2}{*}{$\hat{\theta}_1 < \hat{\theta}_2 < \hat{\theta}_3$} 
    & $\hat{\theta}_2>1$ 
    & 0 
    & 0 \\
    (2) &
    & $\hat{\theta}_2 < 1$ 
    & $\left(1 - \frac{p_{n}}{(1+\beta -\alpha \beta)\delta V}\right) \cdot \lambda$ 
    & 0 \\
    \hline 
    (3) &
    \multirow{3}{*}{$\hat{\theta}_3 < \hat{\theta}_2 < \hat{\theta}_1$} 
    & $1 < \hat{\theta}_3$ 
    & 0 
    & 0 \\
    (4) &
    & $\hat{\theta}_3< 1 < \hat{\theta}_1$ 
    & 0 
    & $\left(1-\frac{p_{r}}{\alpha \delta V}\right) \cdot \lambda$ \\
    (5) &
    & $\hat{\theta}_1<1$
    & $\left(1-\frac{p_{n}-p_r}{(1+\beta)(1-\alpha)\delta V}\right) \cdot \lambda$ 
    & $\left(\frac{p_{n}-p_r}{(1+\beta)(1-\alpha)\delta V} - \frac{p_{r}}{\alpha \delta V}\right) \cdot \lambda$ \\
    \hline
\end{tabular}
\end{table}

Further, we reformulate the expressions of $\hat{\theta}$s as functions between $p_n$ and $p_r$, and display them in the coordinate system about $p_n$ and $p_r$, see Fig. 3. %~\ref{fig: Model O demand rates cases}. 
We merge cases (1) and (3) because they have the same values of $\Lambda_n$ and $\Lambda_r$. 
As illustrated in this figure, the boundary functions divide the range of $p_r \leq p_n$ into four distinct regions, and each region indicates a market outcome:
\begin{itemize}
\item Region \uppercase\expandafter{\romannumeral1} (infeasible market): Neither product type has market demand.
\item Region \uppercase\expandafter{\romannumeral2} (new-only market): Only new products exhibit market demand.
\item Region \uppercase\expandafter{\romannumeral3}(remanufactured-only market): Only remanufactured products show positive demand.
\item Region \uppercase\expandafter{\romannumeral4}(co-existence market): Both new and remanufactured products have market demand.
\end{itemize}

\begin{figure}
\centering
\includegraphics[width = 0.7\linewidth]{images/region chart by price.pdf}
\caption{Region chart of $\Lambda_n$ and $\lambda_r$ with $\beta^-\in [-1,0]$} \label{fig: Model O regions appendix}
\end{figure}

In Fig. 3, %~\ref{fig: Model O demand rates cases},
we reformulate the boundary functions as the relationships between consumer utilities, in order to be consistent with the focus of the previous analysis.

\subsection{Proposition 1}
% ~\ref{pro: Model O approximate optimal prices}}
\label{app: Model O approximate optimal prices}
To achieve the approximation of optimal solutions of Model O, we discuss all cases of $\Lambda_{n}$ and $\Lambda_r$. We first maximize the approximate expected profit for each region shown in Fig.~\ref{fig: Model O regions appendix} and subsequently determine the overall maximum profit. 

We start by Region \uppercase\expandafter{\romannumeral4}, where $\Lambda_{n} = (1-\frac{p_{n}-p_r}{(1+\beta)(1-\alpha)\delta V})\cdot \lambda$ and $\Lambda_r = (\frac{p_{n}-p_r}{(1+\beta)(1-\alpha)\delta V} - \frac{p_{r}}{\alpha \delta V}) \cdot \lambda$. 
By substituting these expressions into the approximate function $\widetilde{\Pi}^O(p_{n},p_r) = (p_{n}-c) \cdot \Lambda_{n}+(p_{r}-c_r-c_{coll})\cdot \Lambda_{r}$, we obtain:
\begin{align*}
    \widetilde{\Pi}^O(p_n,p_r) &= (p_{n} - c) \cdot \left(1 - \frac{p_{n} - p_r}{(1+\beta)(1-\alpha)\delta V }\right) \lambda + (p_r - c_r - c_{coll}) \cdot \left(\frac{p_{n} - p_r}{(1+\beta)(1-\alpha)\delta V} - \frac{p_r}{\alpha \delta V}\right)\lambda \\
    &=  - c \lambda + \frac{\lambda}{\alpha (1+\beta) (1-\alpha)\delta V} 
    \bigg[
    -\alpha \cdot p^2_{n} - (1+\beta-\alpha \beta) \cdot p_r^2 + 2 \alpha \cdot p_{n} p_r \\
    & \quad \quad + \alpha \cdot \Big( c - c_r - c_{coll} + (1+\beta)(1-\alpha) \delta V \Big) \cdot p_{n}
    - \Big( \alpha c - (c_r+c_{coll})(1+\beta-\alpha\beta)\Big) \cdot p_r \notag \bigg].
\end{align*}

The first- and second- order derivatives of $\widetilde{\Pi}^O(p_n,p_r) $ with respect to $p_n$, $p_r$ are as follows:
\begin{align*}
&\frac{\partial \widetilde{\Pi}^O(p_n,p_r) }{\partial p_{n}} = \frac{- \lambda}{(1+\beta)(1-\alpha)\delta V} \bigg[2p_n - 2p_r - (c - c_r - c_{coll}) - (1+\beta)(1-\alpha)\delta V \bigg],\\
&\frac{\partial \widetilde{\Pi}^O(p_n,p_r) }{\partial p_r} = \frac{-\lambda}{\alpha (1+\beta)(1-\alpha)\delta V} \bigg[2(1+\beta-\alpha\beta) \cdot p_r - 2\alpha \cdot p_{n} + \alpha c - (1+\beta-\alpha\beta)(c_r+c_{coll})\bigg], \\
&\frac{\partial^2 \widetilde{\Pi}^O(p_n,p_r) }{\partial p_{n} ^2} = \frac{-2 \lambda}{(1+\beta)(1-\alpha)\delta V}, \notag \\
&\frac{\partial^2 \widetilde{\Pi}^O(p_n,p_r) }{\partial p_r ^2} = \frac{-2 \lambda}{\alpha(1+\beta)(1-\alpha)\delta V} (1+\beta-\alpha\beta). \notag 
\end{align*}

Then we derive the Hessian matrix is as follows:
\begin{equation*}
H=\left[\begin{array}{cc}
\frac{\partial^2 \widetilde{\Pi}^O}{\partial p_n^2} & \frac{\partial^2 \widetilde{\Pi}^O}{\partial p_n \partial p_r} \\
\frac{\partial^2 \widetilde{\Pi}^O}{\partial p_r \partial p_n} & \frac{\partial^2 \widetilde{\Pi}^O}{\partial p_r^2}
\end{array}\right]= \frac{2\lambda}{\alpha(1+\beta)(1-\alpha)\delta V} \left[\begin{array}{cc}
-\alpha  & \alpha \\
\alpha & -(1+\beta-\alpha\beta) \\
\end{array}\right].
\end{equation*}

Because $\frac{-2\lambda}{(1+\beta)(1-\alpha)\delta V}<0$ and $\frac{4\lambda^2 t^2}{\alpha (1+\beta)(1-\alpha)\delta^2 V^2}>0$, 
the Hessian matrix is negative-definite and $\widetilde{\Pi}^O(p_n,p_r)$ is concave. Solving the first conditions, that is $\frac{\partial \widetilde{\Pi}^O(p_n,p_r) }{\partial p_{n}}=0,\frac{\partial \widetilde{\Pi}^O(p_n,p_r) }{\partial p_r} =0$, we obtain the critical point and corresponding maximized approximate profit of Region \uppercase\expandafter{\romannumeral4}: 
\begin{align}
&\widetilde{p_n}^{O} = \frac{c+(1+\beta-\alpha\beta)\delta V}{2}, \notag \\
&\widetilde{p_r}^{O} = \frac{c_r+c_{coll}+\alpha\delta V}{2}, \notag \\
&\widetilde{\Pi}^{O*} = 
    \frac{\lambda}{4 \delta V} \Big[
    \frac{\big((1+\beta^-) (1-\alpha) \delta V-(c+c_r+c_{coll})\big)^2}{(1+\beta^-) (1-\alpha)} + \frac{\big( \alpha\delta V - c_r-c_{coll}\big)^2}{\alpha} \Big]. \label{eq:max profit region 4}
\end{align}

As there is only one product on the market in regions \uppercase\expandafter{\romannumeral2} and \uppercase\expandafter{\romannumeral3}, we could solve the critical points in these Regions based on the analog proof with Model N. The region \uppercase\expandafter{\romannumeral1} has zero demand for new and remanufactured products, so it leads to zero profit. We
summarize the critical points of each region in Table~\ref{tab:Approximate prices and profits}.

\begin{table}\small
\centering
\caption{Critical points and maximized approximate profits}\label{tab:Approximate prices and profits}
\begin{tabular}{|c|c|c|c|}
    \hline 
    \textbf{Regions} & \textbf{$\widetilde{p_n}^{O}$} & \textbf{ $\widetilde{p_r}^{O}$} & \textbf{$\widetilde{\Pi}^{O}$}\\
    \hline 
    \uppercase\expandafter{\romannumeral1} 
    & $\infty$ 
    & $\infty$ 
    & $0$ \\
    \uppercase\expandafter{\romannumeral2} 
    & $\frac{1}{2}\Big(c+(1+ \beta - \alpha \beta)\delta V\Big)$
    & $\infty$ 
    & $\frac{\lambda}{4 (1+\beta -\alpha\beta)\delta V}(c - (1+\beta -\alpha\beta)\delta V)^2$ \\
    \uppercase\expandafter{\romannumeral3} 
    & $\infty$ 
    & $\frac{1}{2}\big(c_r+c_{coll}+\alpha \delta V\big)$
    & $\frac{\lambda}{4 \alpha \delta V}(c_r + c_{coll} - \alpha \delta V)^2$ \\
    \uppercase\expandafter{\romannumeral4}
    & $\frac{1}{2}\Big(c+(1+\beta-\alpha\beta)\delta V \Big)$
    & $\frac{1}{2}\Big(c_r+c_{coll}+\alpha\delta V \Big)$
    & Eq.~\eqref{eq:max profit region 4} \\
    \hline
    \end{tabular}
\end{table}

% Given the critical points (i.e., local maxima) of each region, we verify if they are the global maximum points of their region, and find the overall maximum point of the range of $p_r \leq p_n$. 

We find Region \uppercase\expandafter{\romannumeral4} is important because its local maximized profit by Eq.~\eqref{eq:max profit region 4} is the sum of that in Region \uppercase\expandafter{\romannumeral2} and \uppercase\expandafter{\romannumeral3}. It indicates that if Region \uppercase\expandafter{\romannumeral4}'s local maximum represents the global maximum of this region, it naturally provides higher profit than other regions and thus becomes the overall maximum. Otherwise, the Region \uppercase\expandafter{\romannumeral4}'s global maximum must occur at one of the boundaries, leading to either new-only market (same as Region \uppercase\expandafter{\romannumeral2}) or remanufactured-only market (same as Region \uppercase\expandafter{\romannumeral3}).

To verify if Region \uppercase\expandafter{\romannumeral4}'s local maximum is its global maximum, we substitute the local maximum $(\frac{c+(1+\beta-\alpha\beta)\delta V}{2}, \frac{c_r+c_{coll}+\alpha\delta V}{2})$ into the boundary conditions, which are $p_r>p_n-(1+\beta)(1-\alpha)\delta V$ and $p_r<\frac{\alpha}{1+\beta-\alpha \beta}p_n$. By simplifying the expressions, we get $\frac{(1+\beta)(1-\alpha)}{\alpha}(c_r+c_{coll})<c-c_r-c_{coll}<(1+\beta)(1-\alpha)\delta V$.  It indicates that if the cost conditions are satisfied, Region \uppercase\expandafter{\romannumeral4}'s local maximum represents the global maximum of the region, and the overall maximum of the whole range because of producing larger profit.

However, if the cost conditions are not satisfied,  two scenarios emerge: 
\begin{itemize}
    \item Scenario 1 (Lower Boundary Violation): If $c-c_r-c_{coll} \leq \frac{(1+\beta)(1-\alpha)}{\alpha}(c_r+c_{coll})$, the local maximum violates the lower boundary due to $\widetilde{p_r}^{O*}\geq \frac{\alpha}{1+\beta-\alpha \beta}\widetilde{p_n}^{O*}$. This means Region \uppercase\expandafter{\romannumeral4}'s global maximum moves to its lower boundary, which leads to new-only market, and is indifferent from Model N.
    \item Scenario 2 (Upper Boundary Violation): If $c-c_r-c_{coll}\geq (1+\beta)(1-\alpha)\delta V$, the local maximum violates the upper boundary due to $\widetilde{p_r}^{O*}\leq \widetilde{p_n}^{O*}-(1-\alpha)(1+\beta)\delta V$. This means Region \uppercase\expandafter{\romannumeral4}'s global maximum moves to its upper boundary, which leads to remanufactured-only market. Therefore, the overall maximum approximate expected profit is $\widetilde{\Pi}^{O*} = \frac{\lambda}{4 \alpha \delta V}  \left( c_r+c_{coll}-\alpha \delta V\right)^2$ and the approximate price of remanufactured product is $\widetilde{p_r}^{O*} = \frac{1}{2} \left(c_r+c_{coll}+\alpha\delta V \right)$.
    % , which could happen as $\alpha$ increases to 1
\end{itemize}

\subsection{Proposition 2}
% ~\ref{pro:selection map model N vs O}}
\label{app:model n vs o}

We compare the approximate maximized profits between Model N and Model O to determine the most profitable model. 

Recall the approximate maximized profit of Model N is $\widetilde{\Pi}^{N*}(\cdot) = \frac{(c-\delta V)^2}{4 \delta V}\lambda$. According to Proposition 2,
% ~\ref{pro: Model O approximate optimal prices},
Model O has two profit expressions under different cost conditions.

When $\frac{(1+\beta^-)(1-\alpha)}{\alpha}(c_r+c_{coll})<c-c_r-c_{coll}< (1+\beta^-)(1-\alpha)\delta V$, 
\begin{equation*}
    \widetilde{\Pi}^{O*} = \frac{\lambda}{4 \delta V} \Big[ \frac{\big((1+\beta^-) (1-\alpha) \delta V-(c+c_r+c_{coll})\big)^2}{(1+\beta^-) (1-\alpha)} + \frac{\big( \alpha\delta V - c_r-c_{coll}\big)^2}{\alpha} \Big].
\end{equation*}

When $c-c_r-c_{coll}\geq (1+\beta^-)(1-\alpha)\delta V$, 
\begin{equation*}
    \widetilde{\Pi}^{O*} = \frac{\lambda}{4 \alpha \delta V}  \big(c_r+c_{coll}-\alpha \delta V\big)^2.
\end{equation*}
Here, we exclude the case where $c-c_r-c_{coll}\leq \frac{(1+\beta)(1-\alpha)}{\alpha}(c_r+c_{coll})$ because its outcome is equivalent to Model N. 

Now we solve $\widetilde{\Pi}^{O*} \geq \widetilde{\Pi}^{N*}$ for each case. By simplifying inequalities, we deduce three distinct ranges based on $\alpha$:
\begin{itemize}
    \item When $0\leq \alpha< \frac{c_r+c_{coll}}{c}$,  $\widetilde{\Pi}^{O*} \geq \widetilde{\Pi}^{N*}$ has no solution for any $|\beta^-| \in [0,1]$.
    \item When $\frac{c_r+c_{coll}}{c} \leq \alpha< \frac{(\delta V-c)^2 + (\delta V-c)\sqrt{(\delta V-c)^2+4\delta V (c_r+c_{coll})}+2\delta V (c_r+c_{coll})}{2 (\delta V)^2}$, $\widetilde{\Pi}^{O*} \geq \widetilde{\Pi}^{N*}$ exists when $|\beta^-| \in [0,|\frac{\alpha c^2-(c_r+c_{coll})^2-\alpha(1-\alpha)(\delta V)^2+\sqrt{\left( (\alpha \delta V+c_r+c_{coll})^2-\alpha(\delta V+c)^2 \right) \cdot \left( (\alpha \delta V-c_r-c_{coll})^2-\alpha(\delta V-c)^2 \right)}}{2\alpha(1-\alpha)(\delta V)^2}|]$.
    \item When $\frac{(\delta V-c)^2 + (\delta V-c)\sqrt{(\delta V-c)^2+4\delta V (c_r+c_{coll})}+2\delta V (c_r+c_{coll})}{2 (\delta V)^2} \leq \alpha \leq 1$, $\widetilde{\Pi}^{O*} \geq \widetilde{\Pi}^{N*}$ holds for any $|\beta^-|\in[0,1]$.
\end{itemize}

To simplify the notations, we define the thresholds in Table~\ref{tab:thresholds}.
\begin{table}[ht!]\small
\centering
\caption{Analytical expression for thresholds in Proposition 2}
% ~\ref{pro:selection map model N vs O}} 
\label{tab:thresholds}
\begin{tabular}{ll}
    \toprule
    \textbf{Notation} & \textbf{Expressions} \\
    \midrule
    $\alpha_1$ & $\frac{c_r+c_{coll}}{c}$ \\
    $\alpha_2$ & $\frac{(\delta V-c)^2 + (\delta V-c)\sqrt{(\delta V-c)^2+4\delta V (c_r+c_{coll})}+2\delta V (c_r+c_{coll})}{2 (\delta V)^2}$\\
    $\beta_1(\alpha)$ & $|\frac{\alpha c^2-(c_r+c_{coll})^2-\alpha(1-\alpha)(\delta V)^2+\sqrt{\left( (\alpha \delta V+c_r+c_{coll})^2-\alpha(\delta V+c)^2 \right) \cdot \left( (\alpha \delta V-c_r-c_{coll})^2-\alpha(\delta V-c)^2 \right)}}{2\alpha(1-\alpha)(\delta V)^2}|$\\
    % $\overline{\beta_1}=\beta_1(\alpha_2)$ & $|\frac{(\delta V-c) \left(\sqrt{(\delta V)^2-2 \delta V c+4 \delta V (c_{coll}+c_r)+c^2}-3 \delta V+c\right)}{2\delta V (2 \delta V-c-c_{coll}-c_r)}|$
    \bottomrule
    \end{tabular}{}
\end{table}

%%%%%%%%%%%%%%%%%%%%%%%%%%%%%%%%%
% Numerical Results
\section{Numerical Results} \label{app: numerical}
We use the following parameter settings in numerical experiments: $\lambda =1000$, $V=1000$, $\delta=0.8$, $c,c_{coll},c_r=200,40,80$, $H =10000$, $h=100$, $\alpha\in\{0,0.01,0.02,\cdots,1\}$, $\beta^-\in\{0,-0.01,-0.02,\cdots,-1\}$, $\beta^+\in\{0,0.01,0.02,\cdots,1\}$.

Table~\ref{tab: approximation gap} shows the relative error between the exact and approximate solutions. Table~\ref{tab:thresholds values} shows the numerical values of the threshold in Fig. 6. %~\ref{fig:selection map H10k h100}.
\begin{table}[ht!]\small
\caption{Approximation error of model N and O} \label{tab: approximation gap}
\centering
\begin{tabular}{@{}lllll@{}}
    \toprule 
    &  &\textbf{Exact} &\textbf{Approximation} &\textbf{\small Relative Error} \\ 
    \midrule
    \multirow{3}{*}{\textbf{Model N}} 
    & Optimal price & 497.74 & 500 & $\approx$ 0.45\% \\
    & Optimal quantity & 383 & 380 & $\approx$ 0.78\% \\
    & Maximized profit & 112488.44 & 112500 & $\approx$ 0.01\% \\
    \midrule
    \multirow{3}{*}{\textbf{Model O}} 
    & Optimal prices     & (492.3,380.0) & (492.0,380.0) & \textless{}0.061\% \\
    & Optimal quantities & (224,193) &(226,190) & \textless{}0.9\% \\
    & Maximized profit & 112692.76 &112736.11 & \textless{}0.04\% \\
    \bottomrule
\end{tabular}{}
\end{table}

\begin{table}[ht!]\small
\centering
\caption{Values of thresholds in Fig. 6} %Fig.~\ref{fig:selection map H10k h100}} 
\label{tab:thresholds values}
\begin{tabular}{ll}
\toprule    \textbf{Notation} & \textbf{Numerical values}\\
\midrule
$\alpha_1$ & $0.6$ \\
    $\alpha_2$ & $0.836$\\
    $\beta_1(\alpha_2)$ & $0.392$\\
    $\alpha_1'$ & $0.45$ \\
    \bottomrule
\end{tabular}
\end{table}